\documentclass[onefignum,onetabnum]{siamart171218}

\usepackage[utf8]{inputenc}
\usepackage[margin=1in]{geometry}

\usepackage{mathtools} %
\usepackage{amssymb}
\usepackage{physics}
\usepackage{bm}
\usepackage{makecell}
\usepackage{graphicx}
\usepackage{subcaption}
\usepackage{array}
\usepackage{breqn}
\usepackage{dsfont}
\usepackage{booktabs}
\usepackage{setspace}

\usepackage{xcolor}
\usepackage{pgfplots}
\usepgfplotslibrary{groupplots,dateplot}
\pgfplotsset{compat=newest}
\usepackage{tikz}
\usetikzlibrary{patterns,shapes.arrows,external,math}

\definecolor{crimson2152528}{RGB}{215,25,28}
\definecolor{darkgray176}{RGB}{176,176,176}
\definecolor{lightgray204}{RGB}{204,204,204}
\definecolor{palegreen171221164}{RGB}{171,221,164}
\definecolor{sandybrown25317497}{RGB}{253,174,97}
\definecolor{steelblue43131186}{RGB}{43,131,186}

\DeclareMathOperator*{\argmax}{arg\,max}

\newcommand{\nmax}{n_{\mathrm{max}}}
\newcommand{\tmax}{t_{\mathrm{max}}}
\newcommand{\jmax}{j_{\mathrm{max}}}
\newcommand{\tmin}{t_{\mathrm{min}}}
\newcommand{\Tmin}{T_{\mathrm{min}}}
\newcommand{\Rmax}{R_{\mathrm{max}}}
\newcommand{\vbnres}{\vb{n}_{\mathrm{res}}}
\newcommand{\bessint}{\mathcal{I}}
\newcommand{\dbstext}[1]{\ \ \text{#1}\ \ }
\newcommand{\ex}[1]{\times 10^{#1}}
\newcommand{\sigmax}{\sigma_{\mathrm{max}}}

\newcommand{\ReviewerOne}[1]{#1}
\newcommand{\ReviewerTwo}[1]{#1}
\newcommand{\ReviewerOneStrike}[1]{}
\newcommand{\ReviewerTwoStrike}[1]{}

\title{Lattice Green's Functions for high order finite difference stencils}
\author{James Gabbard and Wim M. van Rees}
\date{January 2023}

\newsiamremark{remark}{Remark}
\newsiamremark{hypothesis}{Hypothesis}
\crefname{hypothesis}{Hypothesis}{Hypotheses}
\newsiamthm{claim}{Claim}

\headers{High order Lattice Green's Functions}{James Gabbard and Wim M. van Rees}

\title{Lattice Green's Functions for high order finite difference stencils\thanks{\textbf{Funding:}~This work was funded by Department of Energy Advanced Scientific Computing Research Program, award no.~DE-SC0020998}}

\author{James Gabbard\thanks{Department of Mechanical Engineering, Massachusetts Institute of Technology (\email{jgabbard@mit.edu}).}
\and Wim M. van Rees\thanks{Department of Mechanical Engineering, Massachusetts Institute of Technology (\email{wvanrees@mit.edu}).}
}

\usepackage{amsopn}

\makeatletter
\newcommand*{\addFileDependency}[1]{%
  \typeout{(#1)}%
  \@addtofilelist{#1}%
  \IfFileExists{#1}{}{\typeout{No file #1.}}%
}
\makeatother

\usepgfplotslibrary{external}
\ifpdf
\hypersetup{
  pdftitle={Lattice Green's Functions for high order finite difference stencils},
  pdfauthor={J. Gabbard and W. M. van Rees}
}
\fi

\begin{document}

\maketitle

\begin{abstract}
    Lattice Green's Functions (LGFs) are fundamental solutions to discretized linear operators, \ReviewerTwo{and as such they are a useful tool for solving discretized elliptic PDEs on domains that are unbounded in one or more directions. The majority of existing numerical solvers that make use of LGFs rely on a second-order discretization and operate on domains with free-space boundary conditions in all directions. Under these conditions, fast expansion methods are available that enable precomputation of 2D or 3D LGFs in linear time, avoiding the need for brute-force multi-dimensional quadrature of numerically unstable integrals.} Here we focus on higher-order discretizations of the Laplace operator on domains with more general boundary conditions, by (1) providing an algorithm for fast and accurate evaluation of the LGFs associated with high-order dimension-split centered finite differences on unbounded domains, and (2) deriving closed-form expressions for the LGFs associated with both dimension-split and Mehrstellen discretizations on domains with one unbounded dimension. Through numerical experiments we demonstrate that these techniques provide LGF evaluations with near machine-precision accuracy, and that the resulting LGFs allow for numerically consistent solutions to high-order discretizations of the Poisson's equation on fully or partially unbounded 3D domains.

\end{abstract}



\section{Introduction}\label{sec:introduction}

A Lattice Green's function (LGF) is a fundamental solution to a discretized linear operator. LGFs are common kernels in fast multipole methods \ReviewerOne{(FMMs)} for elliptic difference equations \cite{gillman2010fast, liska2014parallel, dorschner2020fast}, and they can serve as a regularized Green's function for \ReviewerOne{Poisson solvers based on the fast Fourier transform (FFT) \cite{chatelain2010fourier, caprace2021flups, balty2023flups}}. These fast elliptic solvers are a key component of incompressible flow simulations on unbounded (or free-space) domains, where the Green's function convolution approach enables the enforcement of far-field boundary conditions outside of the compact solution domain \cite{liska2016fast, beckers2022planar, eldredge2022method}. For these methods the LGF, as opposed to a continuous Green's function, can be used to enforce discrete conservation properties \cite{liska2016fast, gabbard2022immersed}, or to preserve the accuracy of immersed interface methods that are tailored to a particular difference scheme \cite{gillis2018fast, gabbard2022immersed}.

\ReviewerOne{The choice of fast convolution algorithm (FFT or FMM) typically depends on the source distribution and the choice of computational domain. For PDEs with smooth source distributions that can be efficiently contained in a rectangular region, FFT-based convolution algorithms for uniform grids outperform the FMM when implemented on massively parallel computer architectures \cite{gholami2016fft}. When the distribution of sources is sparse or irregularly shaped, the flexibility of the FMM can become advantageous, and there have been several FMMs that make use of LGFs to solve elliptic difference equations using uniform grids \cite{gillman2010fast, liska2014parallel} or multiresolution grids \cite{dorschner2020fast}. In either case, the accuracy of  these convolution algorithms is limited by the accuracy of the kernel, so that a solution that satisfies a discretized PDE to near machine precision requires an evaluation of the LGF to near machine precision. Thus the range of difference equations that can be solved with either FFT-based or FMM-based approaches depends on the availability of algorithms that can accurately evaluate the associated LGFs.}

In prior works, the LGF has been synonymous with the LGF for the second-order centered difference stencil for the Laplace operator on a fully unbounded 2D or 3D Cartesian grid. This difference scheme arises also in lattice models with nearest-neighbor interactions, and as a result there has been a considerable amount of work in the physics community on evaluating or approximating specific values of this LGF with analytical techniques; see \cite{zucker201170+} for an \ReviewerTwoStrike{excellent} overview. One such result is that the second-order LGF on a fully unbounded domain can be expressed as a 1D improper integral involving modified Bessel functions \cite{koster1954simplified, cserti2000application}, which has allowed for fast and \ReviewerOneStrike{precise }\ReviewerOne{accurate} evaluations of the second-order LGF in numerical solvers. \ReviewerTwo{However, similar results have not been proposed for the LGFs associated with higher-order finite difference stencils}, which are generally not analogous to common physical models and do not share the same evaluation pathways that rely on Bessel functions.

Additionally, while most existing literature has focused on the LGF for unbounded domains, the concept can be extended to domains with periodic or symmetric boundary conditions along one or more boundaries. These boundary conditions  occur frequently in aerodynamics applications: symmetric-unbounded domains allow for robust inflow and outflow boundary conditions in the streamwise direction while eliminating blockage effects from the other boundaries, \cite{chatelain2013} while periodic-unbounded domains allow for studies of isolated vortex tubes without artificial boundary conditions \cite{ji2022bursting}. \cite{buneman1971analytic} provides a closed-form solution for the second-order stencil in a two dimensional domain with one unbounded and one periodic direction. \ReviewerTwo{However, no corresponding results have been obtained for the LGFs of higher-order finite difference stencils}.

\ReviewerTwo{Besides the LGFs, an alternative strategy for solving elliptic PDEs that is amenable to higher-order convergence and more general domains relies on regularized Green's functions, which are nonsingular approximation to the Green's function of the continuous Laplace operator. A regularization approach based on cell averaging that is compatible with mixed periodic and unbounded domains is developed in \cite{chatelain2010fourier}, which demonstrates second-order convergence. The work of \cite{spietz2018regularization} extends these techniques to the Gaussian regularization proposed in \cite{hejlesen2013high, hejlesen2015high}, allowing for higher-order convergence in Poisson problems on mixed periodic and unbounded domains. While these regularized Green's functions are an effective method for solving elliptic PDEs on unbounded domains, they do not provide an exact solution for elliptic difference equations derived from these PDEs.}

\ReviewerTwo{Extending LGFs to high-order finite difference stencils and more general domain boundary conditions using existing analytical techniques presents several immediate challenges. In fully unbounded domains the LGF is defined by a Fourier integral with a singular integrand. A far-field expansion of this integral for arbitrary finite difference schemes can be computed using the algorithm in \cite{martinsson2002asymptotic}; however, near-field interactions must be computed directly via multidimensional quadrature or approximated with ad-hoc methods \cite{kavouklis2018computation}. For 3D domains with one or two unbounded directions, FFT-based solvers require a partially Fourier-transformed Green's function which can also be expressed as a Fourier integral. Here the integrand is non-singular, and there is no general method that yields a far-field expansion. Evaluating far-field interactions via quadrature leads to catastrophic cancellation due to a rapidly oscillating integrand, while existing workarounds developed for regularized Green's functions lead to only an approximate enforcement of free-space boundary conditions \cite{spietz2018regularization}.}

\ReviewerTwo{In this work we address two major challenges inherent in evaluating the LGFs of high-order finite difference schemes for the Laplace operator. First, we provide an algorithm that can evaluate near-field LGF values for arbitrary dimension-split stencils on fully unbounded domains, yielding near machine-precision results and avoiding the need for multidimensional quadrature. Second, we develop an analytical toolset that yields numerically stable closed-form expressions for the LGF of an arbitrary difference stencil on domains with one unbounded direction. When combined with the far-field expansion algorithms developed in \cite{martinsson2002asymptotic}, our first contribution allows the fast and accurate precomputation of LGFs of higher-order finite difference stencils on fully unbounded domains.  Our second contribution extends this capability to domains with one unbounded direction and one or two periodic directions. The precomputed LGFs can then be paired with either FFT-based or FMM-based convolution algorithms to obtain exact solutions to high-order discretizations of elliptic PDEs.}

\ReviewerTwo{The remainder of this work is organized as follows.} We begin by defining our notation and the high-order finite difference schemes of interest in section \ref{sec:stencils}. Section \ref{sec:unbounded} \ReviewerTwo{presents an algorithm that allows for fast and accurate evaluations of the LGFs associated with dimension-split stencils} on fully unbounded domains. Section \ref{sec:one_unbounded} develops a procedure for deriving stable closed-form expressions for the LGF in domains with a single unbounded dimension, applicable to any consistent and symmetric finite difference discretization. Section \ref{sec:2d} addresses several technicalities that arise when extending these 3D results to 2D domains. We provide computational results that quantify the accuracy and efficiency of these computational methods in section \ref{sec:results}, and draw conclusions in section \ref{sec:conclusion}.

\section{Finite Difference Stencils}\label{sec:stencils}
\ReviewerTwo{In this section we define the notation used throughout this work to represent finite difference operators, both through their coefficients and through their Fourier transforms.} This \ReviewerTwo{work focuses on} finite difference approximations to the negative Laplace operator ($-\nabla^2$) that are defined on a $d$-dimensional Cartesian lattice for $d = 2,3$. We consider lattices with unit spacing, and treat functions on a lattice as functions of an integer index $\vb{n} = (n_1,\,\dots,\,n_d)$. Below we detail the finite-difference schemes considered in this work, starting from a general Mehrstellen discretization in section~\ref{subsec:mehrstellen} and then specializing the discussion to classical dimension-split stencils in section~\ref{subsec:dimsplit}. We note that all finite difference discretizations discussed here are consistent and symmetric, and any bounded function in the null space of the difference operators is a constant function. These requirements place useful restrictions on the spectrum of the difference operator, as discussed in section~\ref{sec:poisson}.

\subsection{Mehrstellen Stencils}\label{subsec:mehrstellen}
A Mehrstellen discretization of the Poisson equation takes the form $\mathcal{L} u = \mathcal{R} f$,  where $u(\vb{n})$ and $f(\vb{n})$ are scalar functions defined on a lattice.
Here $\mathcal{L}$ and $\mathcal{R}$ are a pair of finite difference operators of the form
\begin{equation}
    \qty[\mathcal{L} u](\vb{n}) = \sum_{\bm{\alpha} \in \mathcal{I}(\mathcal{L})} a_{\bm{\alpha}} u (\vb{n} + \bm{\alpha}), \quad \qty[\mathcal{R} f](\vb{n}) = \sum_{\bm{\beta} \in \mathcal{I}(\mathcal{R})} b_{\bm{\beta}} f(\vb{n} + \bm{\beta}),
\end{equation}
where $\bm{\alpha}, \, \bm{\beta} \in \mathbb{Z}^d$ are integer indices, $a_{\bm{\alpha}}, \, b_{\bm{\beta}}$ are the corresponding finite difference coefficients, and $\mathcal{I}(\cdot)$ indicates the set of indices for the coefficients of a finite difference operator. \ReviewerTwo{The design and analysis of Mehrstellen discretizations is treated comprehensively in \cite{deriaz2020compact}, which also provides a brief history of their development. Following \cite{deriaz2020compact},} we define symmetry by requiring that $a_{\bm{\alpha}} = a_{\pi(\bm{\alpha})}$ and $b_{\bm{\beta}} = b_{\pi(\bm{\beta})}$ for any transformation $\pi(\cdot)$ generated by flipping signs or permuting entries within a multi-index. The Fourier symbols of the operators $\mathcal{L}$ and $\mathcal{R}$ are then given by
\begin{equation}
    \sigma_{\mathcal{L}}(\bm{k}) = \sum_{\bm{\alpha} \in \mathcal{I}(\mathcal{L})} a_{\bm{\alpha}} e^{i \bm{k} \cdot \bm{\alpha}}, \quad 
    \sigma_{\mathcal{R}}(\bm{k}) = \sum_{\bm{\beta} \in \mathcal{I}(\mathcal{R})} b_{\bm{\beta}} e^{i \bm{k} \cdot \bm{\beta}}.
\end{equation}
 Consistency implies that $\sigma_{\mathcal{L}}(\bm{k}) = \abs{\bm{k}}^2 + \mathcal{O}(\abs{\bm{k}}^4)$ and $\sigma_{\mathcal{R}}(\bm{k}) = 1 + \mathcal{O}(\abs{\bm{k}}^2)$ as $\abs{\bm{k}} \rightarrow 0$. By taking advantage of symmetry, the symbols can be reduced to the real form
\begin{equation}
\sigma_{\mathcal{L}}(\bm{k}) = \sum_{\bm{\alpha} \in \mathcal{I}(\mathcal{L})} a_{\bm{\alpha}} \prod_{i = 1}^d \cos(\alpha_i k_i),
\end{equation}
and likewise for $\sigma_{\mathcal{R}}(\bm{k})$. While convenient, this form can suffer from a loss of precision as $\abs{\bm{k}} \rightarrow 0$, which can be remedied by re-expressing the symbol as a polynomial in the variables $y_i = \sin^2(k_i/2)$. Using the identities $\cos(\theta) = 1 - 2\sin^2(\theta/2)$ and $\cos(n\theta) = T_n(\cos(\theta))$, where $T_n(x)$ is a Chebyshev polynomial of the first kind, we find the numerically stable form 
\begin{equation}
   \tilde{\sigma}_{\mathcal{L}}(\bm{y}) = \sum_{\bm{\alpha} \in \mathcal{I}(\mathcal{L})} a_{\bm{\alpha}} \prod_{i = 1}^d T_{\abs{\alpha_i}} (1 - 2y_i).
\end{equation}
The same transformation can be applied to the symbol of the right-hand side operator to obtain $\tilde{\sigma}_{\mathcal{R}}(\bm{y})$. Throughout the article will use the Mehrstellen stencils of order four and six \cite{deriaz2020compact}, which have coefficients listed in Table \ref{tab:meh_coefficients}. 

\begin{table}[htb!]
    \centering
    \renewcommand{\arraystretch}{1.0}
    \begin{spacing}{1.3}
    \begin{tabular}{ c m{5.0cm} m{0.1cm} m{5.0cm} }
        \toprule
        Order & $\{a_{\bm{\alpha}}\}$ & & $\{b_{\bm{\beta}}\}$ \\
        \midrule
        4 & $a_{\pi(0,0,0)} = 4$, $a_{\pi(0, 0, 1)} = -\frac{1}{3}$, $a_{\pi(0, 1, 1)} = -\frac{1}{6}$
          & & $b_{\pi(0,0,0)} = \frac{1}{2}$, $b_{\pi(0, 0, 1)} = \frac{1}{12}$ \\
        \addlinespace[0.2em]
        \midrule
        \addlinespace[0.2em]
        6 & {\setlength{\baselineskip}{1.5\baselineskip} $a_{\pi(0,0,0)} = \frac{64}{15}$, $a_{\pi(0, 0, 1)} = -\frac{7}{15}$, $a_{\pi(0, 1, 1)} = -\frac{1}{10}$, $a_{\pi(1, 1, 1)} = -\frac{1}{30}$}
          & & $b_{\pi(0,0,0)} = \frac{67}{120}$, $b_{\pi(0,0,1)} = \frac{1}{18}$, 
          $b_{\pi(0,0,2)} = -\frac{1}{240}$,
          $b_{\pi(0,1,1)} = \frac{1}{90}$,  \\
          \bottomrule
    \end{tabular}
    \end{spacing}
    \caption{Coefficients of the fourth and sixth order Mehrstellen stencils \ReviewerTwo{\cite{deriaz2020compact}}. }
    \label{tab:meh_coefficients}
\end{table}

\subsection{Dimension-split stencils}\label{subsec:dimsplit}
The notation for the Mehrstellen stencils above can be simplified to obtain classical dimension-split discretizations of the Poisson equation, for which $\mathcal{R}$ is the identity and $\mathcal{L}$ is a finite difference operator of the form
\begin{equation}
    \qty[\mathcal{L}u](\vb{n}) = \sum_{i = 1}^d \sum_{j = -w}^w a_j u(\vb{n} + j\vb{e}_i).
\end{equation}
Here $\vb{e}_i$ is a unit vector for the $i$-th dimension, $w$ is the stencil width, and the $\{a_j\}$ are stencil coefficients. Symmetry and consistency imply that $a_{-j} = a_{j}$ and $a_0 = -2 \sum_{j = 1}^{w} a_j$, so that $\mathcal{L}$ is uniquely defined by its non-central coefficients $\{a_j\}_{j = 1}^w$. For these discretizations $\sigma_{\mathcal{R}}(\bm{k}) = 1$ and $\sigma_{\mathcal{L}}(\bm{k}) = \sum_{i = 1}^d \sigma(k_i)$, with
\begin{equation}
\label{eq:symbol_dimension_split}
    \sigma(k) \equiv \sum_{j = -w}^w a_j e^{ijk} = 4 \sum_{j = 1}^w a_j \sin[2](\frac{jk}{2}).
\end{equation}
While the two forms in Equation~\eqref{eq:symbol_dimension_split} are equivalent, the second is more numerically stable for small $k$. By definition $\sigma(k)$ is even and $2\pi$-periodic, and the consistency of the stencil implies that $\sigma(k) = k^2 + \order{k^4}$ at the origin. In addition, we will assume that $\sigma(k)$ is strictly positive for $k \neq 0$, so that any bounded function in the nullspace of $\mathcal{L}$ is a constant function. 

Throughout this work the standard centered stencils of order two, four, six, and eight are used as examples. Their coefficients are given in Table~\ref{tab:split_coefficients}, along with the maximum value of the symbol $\sigma_{\mathrm{max}} = \max_{k} \sigma(k)$. This maximum occurs at $k = \pi$ for all four stencils.
\begin{table}[htb!]
    \centering
    {\renewcommand{\arraystretch}{1.3}%
    \begin{tabular}{c  c  l  c}
        \toprule
        Order & $w$ & $\{a_j\}_{j = 1}^w$ & $\sigma_{\mathrm{max}}$ \\
        \midrule
        2 & 1 & $\qty{-1}$ & $4$ \\
        4 & 2 & $\qty{-\frac{4}{3}, \frac{1}{12}}$ & $\frac{16}{3} \approx 5.33$ \\
        6 & 3 & $\qty{-\frac{3}{2}, \frac{3}{20}, -\frac{1}{90}}$ & $\frac{272}{45} \approx 6.04$ \\
        8 & 4 & $\qty{-\frac{8}{5}, \frac{1}{5}, -\frac{8}{315}, \frac{1}{560}}$ & $\frac{2048}{315} \approx 6.50$ \\
        \bottomrule
    \end{tabular}}
    \caption{Coefficients of the standard second difference stencils. }
    \label{tab:split_coefficients}
\end{table}

\section{Solving the Poisson Equation}\label{sec:poisson}
In this work we consider the discretized 3D Poisson equation $\mathcal{L} u = \mathcal{R} f$ on three types of lattices: those that are fully periodic, those that are fully unbounded, and those that are unbounded \ReviewerTwoStrike{along the first coordinate $n_1$ only} \ReviewerTwo{in only one direction}. In all three cases, the Lattice Green's function (LGF) associated with a finite difference discretization $(\mathcal{L}, \mathcal{R})$ is a scalar function $G(\vb{n})$ satisfying
\begin{equation}\label{eq:lgf_def}
    \qty[\mathcal{L} G](\vb{n}) = \qty[\mathcal{R} \delta](\vb{n}), \qq{with} \delta(\vb{n}) = 
    \begin{cases} 
        1, & \vb{n} = 0 \\ 
        0, & \vb{n} \neq 0.
    \end{cases}
\end{equation}
Once the LGF is known, the solution of the Poisson equation can be obtained from the convolution
\begin{equation}
    u(\vb{n}) = \qty[G \ast f](\vb{n}) \equiv \sum_{\vb{n}'} G(\vb{n} - \vb{n}') f(\vb{n}').
\end{equation}
Here the sum is taken over all lattice points. The LGF is most easily computed via its Fourier transform
\begin{equation}
    \hat{G}(\bm{k}) = \mathcal{F}\qty[G](\bm{k}) \equiv \sum_{\vb{n}} G(\vb{n}) e^{-i\vb{n} \cdot \bm{k}}.
\end{equation}
Applying the transform to both sides of (\ref{eq:lgf_def}) yields the straightforward relation $\hat{G}(\bm{k}) = \sigma_{\mathcal{R}}(\bm{k}) / \sigma_{\mathcal{L}}(\bm{k})$. From here the most efficient algorithm for computing the convolution $u = G \ast f$ varies based on the lattice boundary conditions. As described in \cite{caprace2021flups}, the convolution can be computed dimension-by-dimension, using a fast Fourier transform to perform the convolution along periodic dimensions and a Hockney-Eastwood algorithm \cite{hockney1988computer} to perform the convolution along unbounded dimensions. Thus \ReviewerTwo{fully periodic domains require the evaluation of} $\hat{G}(\bm{k})$; \ReviewerTwo{fully unbounded domains require} the LGF $G(\vb{n})$ obtained via the inverse Fourier transform
\begin{equation}\label{eq:3unb_def}
    G(\vb{n}) = \frac{1}{\qty(2\pi)^3} \int_{[-\pi, \pi]^3} e^{i\vb{n}\cdot\bm{k}} \frac{ \sigma_{\mathcal{R}}(\bm{k}) }{\sigma_{\mathcal{L}}(\bm{k})} \dd{\bm{k}};
\end{equation}
and \ReviewerTwo{domains that are unbounded in only one direction require} the partially transformed Green's function
\begin{equation}\label{eq:1unb_def}
    G(\ReviewerTwo{n}, k_2, k_3) = \frac{1}{2\pi} \int_{-\pi}^{\pi} e^{i \ReviewerTwo{n} k} \frac{ \sigma_{\mathcal{R}}(k, k_2, k_3) }{\sigma_{\mathcal{L}}(k, k_2, k_3)} \dd{k}.
\end{equation}
\ReviewerTwo{Here we have assumed that the first dimension is unbounded, so that $n$ represents an unbounded lattice coordinate and $(k_2, k_3)$ are wavenumbers for the second and third dimensions respectively; other configurations can be obtained by permuting indices.} 

While the computation for periodic domains is trivial, the integrals in Equations~(\ref{eq:3unb_def}) and (\ref{eq:1unb_def}) are challenging to compute numerically. We address Equation~(\ref{eq:3unb_def}) in \ReviewerOne{section \ref{sec:unbounded}}, and return to Equation~(\ref{eq:1unb_def}) in section \ref{sec:one_unbounded}. \ReviewerOne{In each case} we first introduce the general approach for Mehrstellen stencils, and then the specializations and algorithmic optimizations for dimension-split stencils. \ReviewerOne{Table \ref{tab:all_contributions} provides a summary of the existing techniques for evaluating these integrals, as well as the contributions made in this work and outlined in the sections below.}

\begin{table}[htb!]
    \centering
    \renewcommand{\arraystretch}{1.3}
    \ReviewerOne{
    \begin{tabular}{ p{1.3in} p{2.2in} p{2.2in}}
    \toprule
        Stencil Type & Three unbounded directions \newline (near-field interactions) & One unbounded direction \newline (all interactions) \\ 
        \midrule
        Any stencil & $\order{N^3}$ with singular integrand, evaluate (\ref{eq:3unb_def}) via quadrature & $\order{N}$ with singular integrand, evaluate (\ref{eq:1unb_def}) via quadrature \\
        Second-order dimension-split & Reduced to $\order{N}$ via reformulation involving Bessel functions \cite{liska2014parallel} & Reduced to $\order{1}$ via analytical expression \cite{buneman1971analytic} \\
        \midrule
        Mehrstellen stencils & No improvement (section \ref{sec:three_unbounded_mehr}) & Reduced to $\order{1}$ via analytical expression (section \ref{sec:one_unbounded_mehr}) \\ 
        High-order dimension-split & Reduced to $\order{N^2}$ via semianalytical reformulation (section \ref{sec:three_unbounded_split}) & Reduced to $\order{1}$ via analytical expression (section \ref{sec:one_unbounded_split}) \\
        \bottomrule
    \end{tabular}
    }
    \caption{\ReviewerOne{Summary of LGF evaluation strategies and their complexities, including algorithms proposed by previous authors and contributions from this work. The first two rows represent existing techniques, while the remaining rows indicate improved algorithms developed in this work. Complexities are largely illustrative, with $N$ indicating the cost of evaluating a one-dimensional integral with a smooth integrand via quadrature. In practice the quadrature is adaptive, and the exact cost of each integration varies with the integrand, limits of integration, and error tolerance. The far-field interactions on domains with three unbounded directions can be evaluated in $\order{1}$ time for any stencil using existing asymptotic expansion techniques \cite{martinsson2002asymptotic}. }}
    \label{tab:all_contributions}
\end{table}
\subsection{Fully unbounded domains}\label{sec:unbounded}

\subsubsection{Mehrstellen stencils}\label{sec:three_unbounded_mehr}
For a general finite difference discretization $(\mathcal{L}, \mathcal{R})$ on a fully unbounded domain, there is currently no universal strategy for simplifying the 3D integral in (\ref{eq:3unb_def}). 
For moderate values of $\vb{n}$ the integral can be approximated via adaptive numerical quadrature, which is necessary to capture the $1/\abs{\bm{k}}^2$ singularity of the integrand at the origin. For larger values of $\vb{n}$ the integrand becomes highly oscillatory, and evaluation via quadrature becomes impractical. In this regime the LGF can be expanded in a large-$\vb{n}$ asymptotic series using the algorithm provided by Martinsson and Rodin \cite{martinsson2002asymptotic}. In practice the values of $G(\vb{n})$ that require quadrature are precomputed and stored in a lookup table, and the remaining values are computed on-the-fly \ReviewerTwo{from the large-$\vb{n}$ expansion \cite{gillman2010fast}}. While this allows for fast evaluation of $G(\vb{n})$ for any $\vb{n}$, the precomputation can be \ReviewerTwoStrike{quite} \ReviewerTwo{computationally} expensive, particularly when near machine-precision accuracy is required (see section \ref{sec:results}).

\subsubsection{Dimension-split stencils}\label{sec:three_unbounded_split}
When $\mathcal{L}$ is a dimension-split operator and $\mathcal{R}$ is the identity, the cost and accuracy of of this precomputation can be improved \ReviewerTwoStrike{dramatically} \ReviewerTwo{through an additive decomposition of $\sigma_{\mathcal{L}}(\bm{k})$}. Noting that $\sigma_{\mathcal{L}}(\bm{k}) = \sigma(k_1) + \sigma(k_2) + \sigma(k_3)$ is invariant with respect to the transformation $k_i \rightarrow -k_i$, the integral in (\ref{eq:3unb_def}) can be re-expressed in the purely real form
\begin{equation}\label{eq:3unb_integral_real}
    G(\vb{n}) = \frac{1}{\qty(2\pi)^3} \int_{[-\pi, \pi]^3} \frac{\cos(n_1 k_1) \cos(n_2 k_2) \cos(n_3 k_3)}{\sigma(k_1) + \sigma(k_2) + \sigma(k_3)} \dd{k^3}.
\end{equation}
Following \cite{koster1954simplified}, applying the identity $y^{-1} = \int_0^\infty e^{-ty} \dd{t}$ then yields
\begin{equation}\label{eq:no_periodic_product}
    G(\vb{n}) = \int_0^\infty \prod_{i=1}^3 \bessint_{\sigma, n_i}(t) \dd{t},
\end{equation}
where $\bessint_{\sigma, n}(t)$ is a family of special functions associated to the symbol $\sigma(k)$ and an integer index $n$ via 
\begin{equation}\label{eq:def_I_sigma_n_t}
    \bessint_{\sigma, n}(t) \equiv \frac{1}{2\pi} \int_{-\pi}^\pi e^{-t\sigma(k)} \cos(n k) \dd{k}.
\end{equation}
\ReviewerTwo{For convenience, we will often split the integral in  (\ref{eq:no_periodic_product}) into several regions of integration, each notated
\begin{equation}\label{eq:3unb_partial_integration}
    G_{[t_1, t_2]}(\vb{n}) \equiv \int_{t_1}^{t_2} \prod_{i = 1}^3 \bessint_{\sigma, n_i}(t) \dd{t}.
\end{equation}}
When $\mathcal{L}$ represents the standard second-order centered \ReviewerTwo{finite difference} stencil with symbol $\sigma(k) = 2 - 2\cos(k)$, the function $\bessint_{\sigma, n}(t)$is equivalent to $e^{-2t}I_n(2t)$, where $I_n(t)$ is a modified Bessel function of the first kind \ReviewerTwo{\cite{koster1954simplified}}. This reduces the fully 3D integral defining the LGF to
$G(\vb{n}) = \int_0^\infty e^{-6t} I_{n_1}(2t) I_{n_2}(2t) I_{n_3}(2t) \dd{t}$, which is a 1D improper integral with a nonsingular integrand. \ReviewerTwo{To avoid integration over an infinite domain the evaluation is typically split into two regions separated by a large cutoff value $\tmax$, so that $G(\vb{n}) = G_{[0, \tmax]}(\vb{n}) + G_{[\tmax, \infty]}(\vb{n})$. The first term is evaluated directly via 1D numerical quadrature, while the  second is evaluated in closed form using a large-$t$ asymptotic expansion of $I_n(t)$ \cite{liska2014parallel}. This algorithm allows for fast and accurate evaluation of the LGF for the second-order dimension split stencil on a fully unbounded domain, and it has been used in the \ReviewerTwoStrike{vast }majority of numerical applications to date \cite{liska2014parallel, liska2016fast, gillis2018fast, dorschner2020fast, caprace2021flups}.}

For higher-order dimension-split stencils, the modified Bessel function substitution is no longer applicable. Instead we propose a novel algorithm for fast and accurate evaluations of $\bessint_{\sigma, n}(t)$. For moderate values of $t$ the integrand in (\ref{eq:def_I_sigma_n_t}) is smooth and periodic on $[0, 2\pi]$, and the integral can be evaluated to machine precision with numerical quadrature. For large $t$ the integrand is exponentially peaked at the origin, which makes quadrature impractical but allows for the evaluation of $\bessint_{\sigma,n}(t)$ via the large-$t$ asymptotic expansion
\begin{equation}\label{eq:3unb_series_I_b2j}
    \bessint_{\sigma, n}(t) = \frac{1}{\sqrt{4\pi t}} \sum_{j=0}^{J - 1} b_{\sigma, j}(n)t^{-j} + R^\bessint_{\sigma,J}(n, t).
\end{equation}
Here each $b_{\sigma, j}(n)$ is an even polynomial of degree $2j$ with $b_{\sigma,0}(n) = 1$, and $R^\bessint_{\sigma, J}(n, t)$ is the remainder due to truncating the series after $J
$ terms. The full set of coefficients $\{b_{\sigma, j}(n)\}_{j < J}$ can be computed symbolically using the algorithm described in Appendix \ref{section:expand_I_sigma_n_t}. 

\ReviewerTwo{To compute the LGF $G(\vb{n})$ efficiently using this large-$t$ expansion of $\bessint_{\sigma, n}(t)$,} we fix a number of terms $J$ and an evaluation region $[0, \nmax]^3$ in advance, and consider the decomposition 
\begin{equation}\label{eq:3unb_tripartite_G}
    G(\vb{n}) = G_{[0, \tmin]}(\vb{n}) + G_{[\tmin, \Tmin]}(\vb{n}) + G_{[\Tmin, \infty]}(\vb{n}).
\end{equation}
For the first term, both the integral defining $G_{[0, \tmin]}(\vb{n})$ and the integrals defining the $\bessint_{\sigma, n_i}(t)$ are computed via numerical quadrature. For the second term, the integral defining $G_{[\tmin, \Tmin]}(\vb{n})$ is computed via numerical quadrature, while each of the $\bessint_{\sigma, n_i}(t)$ that appear in the integrand are approximated via the large-$t$ series expansion in (\ref{eq:3unb_series_I_b2j}). Finally, the third term is computed in closed form using the expansion
\begin{equation}\label{eq:3unb_outer_expansion}
        G_{[\Tmin, \infty]}(\vb{n}) = \frac{1}{\sqrt{16\pi^3 T}}\sum_{j=0}^{J - 1} g_{\sigma, j}(\vb{n}) \Tmin^{-j} + R^G_{\sigma, J}(\Tmin),
\end{equation}
where $R^G_{\sigma, J}(\Tmin)$ is the remainder due to truncating the series after $J$ terms and each $g_{\sigma, j}(\vb{n})$ is a multivariate polynomial of degree $2j$ defined by
\begin{equation}\label{eq:3unb_outer_gl}
    g_{\sigma,j}(\vb{n}) = \frac{1}{2j+1} \sum_{\ell_1 + \ell_2 + \ell_3 = j} b_{\sigma, \ell_1}(n_1) b_{\sigma, \ell_2}(n_2) b_{\sigma, \ell_3}(n_3).
\end{equation}
\ReviewerTwo{This expansion is obtained directly from (\ref{eq:3unb_partial_integration}) by replacing each $\bessint_{\sigma, n_i}(t)$ with a large-$t$ expansion and integrating each resulting term over the region $[\Tmin, \infty]$.}

\ReviewerTwo{The algorithm is completed by an appropriate choice of the thresholds $\tmin$ and $\Tmin$, which determine the accuracy of the large-$t$ expansion in (\ref{eq:3unb_series_I_b2j}) and (\ref{eq:3unb_outer_expansion}) respectively.} If each $\bessint_{\sigma,n}(t)$ is required with absolute and relative errors $\epsilon_a$ and $\epsilon_r$, we select thresholds $\tmin$ and $\nmax$ such that $\abs{n} \le \nmax$ and $t > \tmin$ implies $\abs{R^\bessint_{\sigma, J}(n, t)} < \epsilon_a$ and $\abs{R^\bessint_{\sigma, J}(n, t)/\bessint_{\sigma, n}(t)} < \epsilon_r$. Assuming that $\abs{\bessint_{\sigma, n}(t)}$ is well represented by the first term in the sum, and that $\abs{R^\bessint_{\sigma, J}(n, t)}$ is well represented by first neglected term, we can estimate $\tmin$ as a function of $J$ and $\nmax$ via
\begin{equation}\label{eq:3unb_lower_threshold}
    \tmin = \max \qty(
        \abs{\frac{b_{\sigma, J}(\nmax)}{\epsilon_r}}^{\frac{1}{J}}, 
        \abs{\frac{b_{\sigma, J}(\nmax)}{\epsilon_a \sqrt{4\pi}}}^{\frac{2}{2J + 1}}
    ).
\end{equation}
\ReviewerTwo{Likewise, given absolute and relative error tolerances $\epsilon_a$ and $\epsilon_r$ for $G_{[\Tmin, \infty]}(\vb{n})$,} we choose a threshold $\Tmin$ such that $T \ge \Tmin$ implies $\abs{R^G_{\sigma, J}(T)} < \epsilon_a$ and $\abs{R^G_{\sigma, J}(T) / G_{[T, \infty]}(\vb{n})} < \epsilon_r$. Assuming that $\abs{G_{[T, \infty]}(\vb{n})}$ is well represented by the first term in the sum and that $\abs{R^G_{\sigma, J}(T)}$ is well represented by the first neglected term, this threshold is approximately
\begin{equation}\label{eq:3unb_upper_threshold}
    T_{\mathrm{min}} = \max \qty( 
        \abs{\frac{g_{\sigma,J}(\vb{n}_{\mathrm{max}})}{\epsilon_r}}^{\frac{1}{J}},
        \abs{\frac{g_{\sigma, j}(\vb{n}_{\mathrm{max}})}{\epsilon_a \sqrt{16\pi^3}}}^{\frac{2}{2J+1}}
    ).
\end{equation}
In the above we choose $\vb{n}_{\mathrm{max}} = [\nmax, \nmax, \nmax]$ under the assumption that the polynomial $g_{\sigma,J}(\vb{n})$ is monotonic in each argument for $n$ above some minimum value. \ReviewerTwo{In our numerical experiments we select a single set of absolute and relative tolerances $\epsilon_a$ and $\epsilon_r$ that are used both for 1D adaptive quadrature routines and for the threshold estimates in ($\ref{eq:3unb_lower_threshold}$) and ($\ref{eq:3unb_lower_threshold}$). Because of the assumption used to derive the thresholds and the additive decomposition in (\ref{eq:3unb_tripartite_G}), the resulting approximation to $G(\vb{n})$ is not strictly bound by the chosen tolerances. However, the numerical experiments presented in section \ref{sec:results} indicate that the resulting errors fall within an order of magnitude of $\epsilon_a$ and $\epsilon_r$.}

\ReviewerTwo{With the algorithm outlined above}, the \ReviewerTwoStrike{vast }majority of the computation is spent evaluating $G_{[0, \tmin]}(\vb{n})$, which as a nested integral is effectively 2D, nonsingular, and has a well-defined domain. \ReviewerOne{The second term $G_{[\tmin, \Tmin]}(\vb{n})$ requires only one-dimensional numerical quadrature, while the third term $G_{[\Tmin, \infty]}(\vb{n})$ is fully analytical.} Compared to a direct evaluation of (\ref{eq:3unb_integral_real}) by 3D adaptive quadrature, the algorithm presented here \ReviewerOne{is able to achieve near machine-precision results with considerably less computation (see section \ref{subsec:residuals_3unb})}.

\subsection{Domains with one unbounded dimension}\label{sec:one_unbounded}
\ReviewerOne{In this section we provide closed-form analytical expressions for LGFs of both Mehrstellen and dimension-split stencils on domains with one unbounded direction, eliminating the need to evaluate (\ref{eq:1unb_def}) via numerical quadrature.} We note that the analysis below is similar in spirit to an approach developed by Buneman \cite{buneman1971analytic} for the second order dimension-split stencil, which relies on the solution of a linear recursion. 
Here we favor an approach based on contour integration, which provides more powerful tools to maintain numerical stability for certain edge cases.

\subsubsection{Mehrstellen stencils}\label{sec:one_unbounded_mehr}
\ReviewerOne{Let $(\mathcal{L}, \mathcal{R})$ be a Mehrstellen discretization on a domain with one unbounded direction.} \ReviewerTwoStrike{Assuming that $n_1$ is the unbounded lattice coordinate, }\ReviewerTwo{Without loss of generality we will assume that $n_1$ is the unbounded lattice coordinate, so that the partially transformed LGF $G(n, k_2, k_3)$ defined in (\ref{eq:1unb_def}) is the object of interest.} We begin by isolating the $k_1$ dependence of $\sigma_{\mathcal{L}}(\bm{k})$. Defining the stencil width $w_{\mathcal{L}} = \max_{\bm{\alpha} \in \mathcal{I}(\mathcal{L})} \norm{\bm{\alpha}}_{\infty}$, the symbol can be rewritten as
\begin{equation}\label{eq:sigma_k_k2_k3}
    \sigma_{\mathcal{L}}(k; k_2, k_3) = \sum_{j = -w_{\mathcal{L}}}^{w_{\mathcal{L}}} a_j(k_2, k_3) e^{ijk}
    \qq{with} a_j(k_1, k_2) = \sum_{\bm{\alpha} \in \mathcal{I}(\mathcal{L}),\,\alpha_1 = j} a_{\bm{\alpha}} e^{i(k_2 \alpha_2 + k_3 \alpha_3)}.
\end{equation}
In this form the symbol resembles that of a one dimensional finite difference scheme with coefficients $\{a_j\}$ that depend on the parameters \ReviewerTwo{$(k_2, k_3)$}. The symmetry of the coefficients $\{a_{\bm{\alpha}}\}$ implies that $a_j(k_2, k_3) = a_{-j}(k_2, k_3)$, so that the one dimensional scheme is symmetric as well. For brevity and numerical stability we replace the parameters $(k_2, k_3)$ with the vector parameter $\bm{y} = [\sin^2(k_2/2),\, \sin^2(k_3/2)] \in [0,\, 1]^2$, so that the coefficients $\{a_j(\bm{y})\}$ are multivariate polynomials in $\bm{y}$. Finally, we define a characteristic polynomial
\begin{equation}\label{eq:1unb_characterstic_poly}
    p_{\mathcal{L}}(z;\bm{y}) = \sum_{j = -w_{\mathcal{L}}}^{w_{\mathcal{L}}} a_j(\bm{y}) z^{j+w_{\mathcal{L}}},
\end{equation}
which has degree $2w_{\mathcal{L}}$ and coefficients equal to the coefficients of the difference scheme. For the discrete operator $\mathcal{R}$, the stencil width $w_{\mathcal{R}}$, coefficients $b_j(\bm{y})$, and characteristic polynomial $p_{\mathcal{R}}(z; \bm{y})$ are defined analogously. 

Returning to the integral in (\ref{eq:1unb_def}), we make the substitution $z = e^{ik}$ and make use of the definitions above to arrive at the contour integral
\begin{equation}\label{eq:1unb_contour}
    G(n; \bm{y}) = \frac{1}{2\pi i} \oint \frac{z^{n - m} p_{\mathcal{R}}(z;\bm{y})}{p_{\mathcal{L}}(z;\bm{y})} \dd{z},
\end{equation}
where $m \equiv w_{\mathcal{R}} - w_{\mathcal{L}} + 1$ and the integral is taken counterclockwise around the unit circle in the complex plane. Taking advantage of the symmetry $G(n;\bm{y}) = G(-n;\bm{y})$, we assume that $n \ge 0$ in (\ref{eq:1unb_contour}) and in the following analysis; for negative arguments each $n$ should be replaced by $\abs{n}$. Because the integrand is a rational function, the integral can be evaluated with the method of residues. Thus the the roots of $p_{\mathcal{L}}(z;\bm{y})$ lying on or inside the unit circle determine the behavior of $G(n, \bm{y})$.

Before proceeding further, we note the assumptions made on the operator $\mathcal{L}$ lead to corresponding constraints on the roots of $p_{\mathcal{L}}(z;\bm{y})$. The consistency of the underlying stencil implies that the polynomial $p(z; \bm{0})$ has a root at $z = 1$ with multiplicity two. To see this we note that $a_{j}(\bm{0}) = \sum_{\alpha_i = j} a_{\bm{\alpha}}$ and that consistency implies that $\mathcal{L}$ is exact on quadratic polynomials, so that
\begin{equation}\label{eq:two_periodic_derivs}
\begin{aligned}
    p_{\mathcal{L}}(1; \bm{0}) &= \sum_{j=-w_{\mathcal{L}}}^{w_{\mathcal{L}}} a_{j}(\bm{0}) = \mathcal{L}[1] = 0, \\ %
    p_{\mathcal{L}}'(1; \bm{0}) &= \sum_{j=-w_{\mathcal{L}}}^{w_{\mathcal{L}}} a_{j}(\bm{0})(j + w_{\mathcal{L}}) = \mathcal{L}[n_1 + w_{\mathcal{L}}] = 0, \\ %
    p_{\mathcal{L}}''(1; \bm{0}) &= \sum_{j=-w_{\mathcal{L}}}^{w_{\mathcal{L}}} a_{j}(\bm{0})(j + w_{\mathcal{L}})(j + w_{\mathcal{L}} - 1) = \mathcal{L}[(n_1 + w_{\mathcal{L}})(n_1 + w_{\mathcal{L}} - 1)] = -2. \\ %
\end{aligned}
\end{equation}
A similar analysis of $\mathcal{R}$ indicates that $p_{\mathcal{R}}(1; \bm{0}) = 1$, which will be useful below. The symmetry of the coefficients of $p_{\mathcal{L}}(z; \bm{y})$ implies that when $r$ is a root of $p_{\mathcal{L}}(z; \bm{y})$, the reciprocal $1/r$ is also a root. Finally, the assumption that $\sigma_{\mathcal{L}}(\bm{k}) \neq 0$ for $\abs{\bm{k}} \neq 0$ implies that $p_{\mathcal{L}}(z; \bm{y})$ is nonzero on the unit circle except when $z = 1$ and $\abs{\bm{y}} = 0$. Taken together, these three properties imply that for $\abs{\bm{y}} > 0$ there are at most $w_{\mathcal{L}}$ roots (counted with their multiplicities) that lie strictly inside the unit circle, and for $\abs{\bm{y}} = 0$ there are at most $w_{\mathcal{L}} - 1$ roots strictly inside the unit circle and a double root at $z = 1$. 

To locate these roots efficiently, let $\lambda = \frac{1}{2}(z + z^{-1})$, and note that $z^j + z^{-j} = 2 T_j(\lambda)$, where $T_j$ is the $j$-th Chebyshev polynomial of the first kind. Making a change of variables in (\ref{eq:1unb_characterstic_poly}) then yields
\begin{equation}\label{eq:1unb_qlambda_definition}
    \frac{p_{\mathcal{L}}(z; \bm{y})}{z^w} = a_0(\bm{y}) + 2 \sum_{j=1}^{w_{\mathcal{L}}} a_j(\bm{y}) T_j(\lambda) \equiv q(\lambda; \bm{y}),
\end{equation}
where $q(\lambda; \bm{y})$ is a degree $w_{\mathcal{L}}$ polynomial in $\lambda$. Differentiating the above definition twice and substituting the relations (\ref{eq:two_periodic_derivs}) yields $q(1; \bm{0}) = 0$ and $q'(1; \bm{0}) = -2$, so that $q(\lambda; \bm{y})$ has a single root at $\lambda = 1$. For each root of $q(\lambda; \bm{y})$, there is a root
\begin{equation}\label{eq:r_from_lambda}
    r = \lambda - \sqrt{\lambda - 1}\sqrt{\lambda + 1}
\end{equation}
satisfying $p_{\mathcal{L}}(r;\bm{y}) = p_{\mathcal{L}}(r^{-1}; \bm{y}) = 0$ with $\abs{r} \le 1$. Thus the change of variables reduces the problem of finding the $2w_{\mathcal{L}}$ roots of $p_{\mathcal{L}}(z;\bm{y})$ to finding the $w_{\mathcal{L}}$ roots of $q(\lambda; \bm{y})$, and allows for a closed form representation of the roots for stencils with $w_{\mathcal{L}} \le 4$. 
Note that (\ref{eq:r_from_lambda}) can be simplified analytically, but the written expression guarantees that $\abs{r} \le 1$ when the square root is taken with a branch cut on the negative real axis.

Returning to the integral in (\ref{eq:1unb_contour}), we begin by assuming that $\abs{\bm{y}} \neq 0$. Let $\{r_i\}$ be the set of roots of $p_{\mathcal{L}}(z;\bm{y})$ lying strictly inside the unit circle. In the simplest case there are $w_{\mathcal{L}}$ distinct roots satisfying $r_i \neq 0$ and  $p_{\mathcal{R}}(r_i; \bm{y}) \neq 0$, so that the integrand has $w_{\mathcal{L}}$ simple poles and a pole of order $m - n$ at $z = 0$. The LGF can then be written in closed form as
\begin{equation}\label{eq:1unb_general}
    G(n; \bm{y}) = \sum_{i = 1}^{w_{\mathcal{L}}} r_i^{n - m} \frac{p_{\mathcal{R}}(r_i; \bm{y})}{p_{\mathcal{L}}'(r_i; \bm{y})} + \dv[m-n-1]{z} \qty[\frac{p_{\mathcal{R}}(z; \bm{y})}{p_{\mathcal{L}}(z; \bm{y})}]_{z = 0},
\end{equation}
where the second is omitted whenever $n > m - 1$. If the set $\{r_i\}$ contains a conjugate pair $(r_1, r_2)$, we use the polar decomposition $r_1 = \rho e^{i\theta}$ to write the corresponding terms in the real form 
\begin{equation}\label{eq:1unb_conjugate_roots}
\begin{aligned}
    \sum_{i = 1}^2 r_i^{n - m} \frac{p_{\mathcal{R}}(r_i; \bm{y})}{p_{\mathcal{L}}'(r_i; \bm{y})} = \abs{K} \rho^n \sin\qty(n\theta  + \phi), \dbstext{with}  K = \frac{2 p_{\mathcal{R}}(r_1; \bm{y})}{r_1^mp_{\mathcal{L}}'(r_1; \bm{y})} \dbstext{and} \phi = \arg \qty(iK).
\end{aligned}
\end{equation}
This form is numerically stable when the roots have positive real part and arbitrarily small imaginary part. 

When $\abs{\bm{y}} = 0$, the integrand in (\ref{eq:1unb_def}) has a $1/k^2$ singularity at the origin, which is not integrable. To remedy this, we instead evaluate
\begin{equation}\label{eq:two_periodic_k}
    G^*(n) = G(n, \bm{0}) - G(0, \bm{0}) = \frac{1}{2\pi} \int_{-\pi}^{\pi} \frac{(e^{ink} - 1)\sigma_{\mathcal{R}}(k, 0, 0)}{\sigma_{\mathcal{L}}(k, 0, 0)} \dd{k}.
\end{equation}
The new integrand is bounded and has a removable singularity at $k = 0$. Expressed as a contour integral,
\begin{equation}\label{eq:two_periodic_Gstar_z}
    G^*(n) = \frac{1}{2\pi i} \oint \frac{(z^n - 1)p_{\mathcal{R}}(z; \bm{0})}{ z^{m} p_{\mathcal{L}}(z; \bm{0})} \dd{z}.
\end{equation}
The denominator has a double root at $z = 1$, but the factorization $z^n - 1 = (z - 1) \sum_{j = 0}^{n - 1} z^j$ indicates that the integrand has only a simple pole at $z = 1$. The corresponding residue is
\begin{equation}\label{eq:1unb_general_residues}
    \lim_{z \rightarrow 1} (z - 1)\qty[ \frac{(z^n - 1)p_{\mathcal{R}}(z; \bm{0})}{ z^{m} p_{\mathcal{L}}(z; \bm{0})} ]
    = \lim_{z \rightarrow 1} \qty[ \frac{p_{\mathcal{R}}(z; \bm{0})}{z^m} \frac{(z - 1)^2}{p_{\mathcal{L}}(z;\bm{0})} \sum_{j = 0}^{n - 1} z^j ]
    = 2n \frac{p_{\mathcal{R}}(1;\bm{0}) }{p_{\mathcal{L}}''(1; \bm{0})} 
    = -n.
\end{equation}
Because this pole lies on the contour of integration, (\ref{eq:two_periodic_Gstar_z}) must be interpreted as the Cauchy principal value
\begin{equation}\label{eq:1unb_y_zero_residues}
    G^*(n) = - \frac{n}{2} + \sum_{i = 1}^{w_{\mathcal{L}} - 1} \frac{(r_i^n - 1)}{r_i^m} \frac{p_{\mathcal{R}}(r_i; \bm{0})}{p_{\mathcal{L}}'(r_i; \bm{0})} + \dv[m-n-1]{z} \qty[\frac{p_{\mathcal{R}}(z; \bm{0})}{p_{\mathcal{L}}(z; \bm{0})}]_{z = 0},
\end{equation}
where the third term is omitted whenever $n > m - 1$. 

We note that for small but nonzero $\bm{y}$ the expression in (\ref{eq:1unb_general_residues}) suffers from catastrophic cancellation in a naive evaluation of $p_{\mathcal{L}}'(z;\bm{y})$. For discretizations of the Poisson equation we do not encounter $\abs{\bm{y}}$ smaller than $\sin^2(\pi/N) \approx \pi^2/N^2$, where $N$ is the number of points along the largest periodic dimension, and result retains a relative precision of roughly $N\epsilon_m$ where $\epsilon_m$ is the relative machine precision. Thus the issue can be safely ignored for most tractable problem sizes. For extremely large problems with high accuracy requirements, Appendix \ref{sec:1unb_small_c} provides a method to restore full relative precision for dimension split stencils with small $\bm{y}$.

Finally, the expressions in (\ref{eq:1unb_general_residues}) and (\ref{eq:1unb_y_zero_residues}) are valid for any $\bm{y} \in [0, 1]^2$ so long as $p_{\mathcal{L}}(z;\vb{y})$ has unique nonzero roots. As these assumptions are relaxed a variety of pathological behavior can occur. In this paper we focus only the Mehrstellen stencils of order four and six, which have $w_{\mathcal{L}} = 1$ and consequently cannot have a repeated root inside the unit circle. However, for the fourth order Mehrstellen stencil the leading coefficient $a_{w_{\mathcal{L}}}(\bm{y})$ vanishes on the set $\mathcal{Y} = \{\bm{y} \in [0,1] \mid y_2 + y_3 = \frac{2}{3}\}$, and the characteristic polynomial $p_{\mathcal{L}}(z; \bm{y})$ has a root at $z = 0$. As $\bm{y}$ approaches the set $\mathcal{Y}$ there is a reciprocal pair of roots $(r,\, r^{-1})$ of $p_{\mathcal{L}}(z;\bm{y})$ that approaches $(0, \infty)$. For $n > m - 1$, the corresponding term in (\ref{eq:1unb_general}) approaches zero, and there is no numerical instability. When $n \le m - 1$ the corresponding term grows without bound, but the $G(n,\bm{y})$ remains finite due to cancellation. To maintain numerical stability, that cancellation must be removed by algebraic manipulation of (\ref{eq:1unb_general}). Noting that $m = 1$ for the fourth order stencil and that both $p_{\mathcal{L}}(z;\bm{y})$ and $p_{\mathcal{R}}(z;\bm{y})$ are quadratic, we make the rearrangement
 \begin{equation}\label{eq:1unb_meh4_stable}
    G(0; \bm{y}) = \qty(\frac{r}{a_1(\bm{y})}) \qty(\frac{b_0(\bm{y}) + 2b_1(\bm{y})r}{r^2 - 1}).
 \end{equation}
 Here we have used the fact that $r$ and $r^{-1}$ are the only roots of $p_{\mathcal{L}}(z;\bm{y})$, so that $p'_{\mathcal{L}}(r; \bm{y}) = a_1(\bm{y}) (r - r^{-1})$. We determine $r$ via the sole root of $q(\lambda; \bm{y})$ given by $\lambda = -a_0(\bm{y})/2 a_1(\bm{y})$; when $\lambda$ becomes large, $r$ can be calculated with the series expansion
\begin{equation}
    r = \frac{1}{2} \lambda^{-1} + \frac{1}{8} \lambda^{-3} + \frac{1}{16} \lambda^{-5} + \frac{5}{128} \lambda^{-7} + \frac{7}{256} \lambda^{-9} + \order{\lambda^{-11}},
\end{equation}
which can be derived from (\ref{eq:r_from_lambda}). In this form 
 the ratio $r/a_1(\bm{y})$ is numerically stable and approaches $-1/a_0(\bm{y})$ as $\bm{y} \rightarrow \mathcal{Y}$, allowing $G(0, \bm{y})$ to be evaluated with full precision.

\subsubsection{Dimension split stencils}\label{sec:one_unbounded_split}
For a dimension-split stencil with coefficients $\{a_j\}_{j = 1}^w$ the polynomial $p_{\mathcal{L}}(z; \bm{y})$ can be expressed in the simplified form
\begin{equation}\label{eq:two_periodic_pzc}
    p_{\mathcal{L}}(z; \bm{y}) = p(z; c) \equiv cz^w + \sum_{j = -w}^{w} a_j z^{j+w},
\end{equation}
where the dependence on $\bm{y}$ is replaced by a dependence on $c = \sigma(k_1) + \sigma(k_2)$, a scalar parameter taking values in the interval $[0, 2\sigmax]$. The polynomial $p(z;c)$ inherits all the properties of $p_{\mathcal{L}}$ described in the previous section. Additionally, because the leading coefficient $a_j$ is constant and nonzero, the polynomial $p(z;c)$ has exactly $w$ nonzero roots (counted with their multiplicities) inside the unit circle for $c \neq 0$, and exactly $w - 1$ roots inside the unit circle for $c = 0$. The corresponding polynomial 
\begin{equation}
    q(\lambda; c) \equiv c + 2\sum_{i = 1}^w a_j \qty(T_j(\lambda) - 1) = \frac{p(z; c)}{z^w}
\end{equation}
depends on $c$ only through the constant term, and can be conveniently written as $q(\lambda) + c$ where $q(\lambda) \equiv q(\lambda; 0)$. When the roots $\{\lambda_i\}$ of $q(\lambda; c)$ are distinct, the roots $\{r_i\}$ of $p(z; c)$ are distinct as well, and after some simplification the LGF can be written in the straightforward form
\begin{equation}\label{eq:1unb_split}
    G(n; c) = \sum_{i = 1}^w \frac{r_i^{n + w - 1}}{p'(r_i; c)} = -\sum_{i = 1}^w \frac{r_i^n}{q'(\lambda_i) \sqrt{\lambda_i - 1}\sqrt{\lambda_i + 1}}.
\end{equation}
 For the dimension-split stencils of order two through eight the polynomials $q(\lambda)$ are listed in Table \ref{tab:qj_lambda}, and the roots of $q(\lambda; c)$ are given explicitly as a function of the parameter $c$ in Appendix \ref{sec:1unb_expressions}.
 
\begin{table}[htb!]
    \centering
    \begin{tabular}{c c}
        \toprule
        Order & $q(\lambda)$ \\
        \midrule
        2 & $-2 \lambda + 2$ \\
        \addlinespace
        4 & $\frac{1}{3}(\lambda^2 - 8\lambda + 7)$ \\ 
        \addlinespace
        6 & $\frac{1}{45}(-4\lambda^3 + 27\lambda^2 - 132\lambda + 109)$ \\
        \addlinespace
        8 & $\frac{1}{315}(9\lambda^4 - 64\lambda^3 + 243\lambda^2 - 960 \lambda + 772)$ \\
        \bottomrule
    \end{tabular}
    \caption{The polynomials $q(\lambda)$  for stencils of order two through eight.}
    \label{tab:qj_lambda}
\end{table}

For the fourth and eighth order dimension split stencils, $p(z; c)$ has a repeated root for a given value of $c > 0$, here denoted $c^*$. The corresponding value $\lambda^*$ is a repeated root of $q(\lambda; c^*)$ as well. For nearby values $c = c^* + \delta$, the repeated root splits into either two real roots $r_1, \, r_2 = \bar{r} \pm \epsilon$ or a conjugate pair $r_1, \, r_2 = \bar{r} \pm i\epsilon$ with $\epsilon \sim \sqrt{\delta}$. This introduces a catastrophic cancellation between the corresponding terms in (\ref{eq:1unb_split}), destroying the precision of the LGF.
For moderately small $\delta$, the terms in (\ref{eq:1unb_split}) corresponding to $r_1$ and $r_2$ can be rewritten to restore numerical stability by taking advantage of the factorization $p(z;c) = a_w \prod_{j = 1}^{2w} (z - r_j)$. After some simplification,
\begin{equation}\label{eq:1unb_repeat_fz}
    \sum_{j = 1}^2 \frac{r_j^{n + w - 1}}{p'(r_j; c)} = \frac{f(r_1) - f(r_2)}{r_1 - r_2} \dbstext{with} f(z) = \frac{z^{n + w - 1}}{a_w \prod_{j = 3}^{2w}(z - r_j)}.
\end{equation}
Considering that $f(z) \sim \mathcal{O}(z^n)$ and $f'(z) \sim \mathcal{O}(n z^n)$ for larger $n$, the expression above has relative precision that is $\mathcal{O}(\epsilon_m / n \sqrt{\delta})$ where $\epsilon_m$ represents the relative machine precision. In double precision arithmetic, this equates to a relative accuracy of $10^{-14}$ whenever $\delta n^2 \gtrsim 10^{-4}$.

Alternatively, the LGF can be expanded in a Taylor series about $c = c^*$. Differentiating (\ref{eq:1unb_split}) with respect to $c$ yields
\begin{equation}\label{eq:1unb_repeated_residues}
    G^{(j)}(n, c^*) = (-1)^j \frac{j!}{2\pi} \int_{-\pi}^{\pi} \frac{\cos(nk)}{ \qty(\sigma(k) + c^*)^{j + 1}} \dd{k} = \frac{(-1)^j j!}{2\pi i} \oint \frac{z^{n + (j + 1)w - 1}}{p(z; c^*)^{j+1}} \dd{z},
\end{equation}
where the superscript in $G^{(j)}$ implies differentiation with respect to $c$. The poles of this integrand are well separated: where $p(z; c^*)$ has a single root, the integrand has a pole of order $j + 1$, and where $p(z;c^*)$ has a repeated root the integrand $I^{(j)}(z, n)$ has a pole of order $2j + 2$. Consequently, (\ref{eq:1unb_repeated_residues}) can evaluated without catastrophic cancellation using the method of residues. For large stencils and large $j$ the computation of residues is tedious and most easily evaluated with a computer algebra system. Once they are known, full LGF can be accurately evaluated with the Taylor expansion
\begin{equation}\label{eq:1unb_repeated_taylor}
    G(n,c) = \sum_{j = 0}^{m - 1} \frac{\delta^j}{j!} G^{(j)}(n, c^*) + \order{\delta^m \bar{r}^n n^{2m+1}}.
\end{equation}
The derivatives $G^{(j)}(n, c^*)$ needed to evaluate (\ref{eq:1unb_repeated_taylor}) for the fourth and eight order dimension split stencils are given explicitly in Appendix \ref{sec:1unb_expressions}.

\subsection{Alterations for 2D domains}\label{sec:2d}
For 2D domains with one unbounded dimension the calculation of the LGF for general or dimension-split stencils is fully analogous to the 3D case. For a fully unbounded 2D domain, the LGF is defined by the integral
\begin{equation}
    G(n_1, n_2) = \frac{1}{(2\pi)^2} \int_{[-\pi,\pi]^2} \frac{e^{i(n_1k_1 + n_2k_2)}}{\sigma(k_1) + \sigma(k_2)} \dd{k_1} \dd{k_2}.
\end{equation}
Here the integrand has a singularity proportional to $1/\abs{k}^2$ at the origin, which is not integrable. The standard strategy for avoiding this singularity is to calculate the relative quantity $G(n_1, n_2) - G(0, 0)$, which is finite for all $n_1$ and $n_2$ \cite{martinsson2002asymptotic}. For a general stencil, the quantity $(e^{i\vb{n}\cdot\bm{k}} - 1)$ replaces $e^{i\vb{n}\cdot\bm{k}}$ in the numerator of (\ref{eq:3unb_def}). For dimension-split stencils the calculation of $\bessint_{\sigma, n}(t)$ remains unchanged, and the LGF is defined by the integral
\begin{equation}
    G(\vb{n}) = \int_0^T \qty[\bessint_{\sigma,n_1}(t) \bessint_{\sigma,n_2}(t) - \bessint_{\sigma,0}(t)^2] \dd{t} + \frac{1}{4\pi} \sum_{j=1}^\infty g_{\sigma, j}(\vb{n}) T^{-j},
\end{equation}
where the multivariate polynomials $g_{\sigma, j}(\vb{n})$ are defined by
\begin{equation}
    g_{\sigma, j}(n_1, n_2) = \frac{1}{j} \sum_{\ell_1 + \ell_2 = j} \qty[b_{\sigma, \ell_1}(n_1) b_{\sigma, \ell_2}(n_2) - b_{\sigma, \ell_1}(0) b_{\sigma, \ell_2}(0)].
\end{equation}
Finally, for 2D domains the far-field expansion provided by Martinsson and Rodin in \cite{martinsson2002asymptotic} is known only up to a fixed constant which varies from stencil to stencil. The constant can be determined by insisting that the evaluation by quadrature and the evaluation by series expansion coincide at some pre-determined point $(n_1, n_2)$ for which both strategies are accurate. 

\section{Results}\label{sec:results}
The computation of LGF values for fully unbounded domains is implemented in Julia within the package \texttt{ExpandLGF.jl}, which is open source and available online\footnote{See github.com/vanreeslab/ExpandLGF.jl}. The calculation of asymptotic expansions relies on the fast symbolic computations provided \texttt{AbstractAlgebra.jl} \cite{AbstractAlgebra.jl-2017}, which are translated into optimized callable Julia functions using \texttt{Symbolics.jl} \cite{gowda2022}. The 1D numerical quadrature used to evaluate the $I_{\sigma, n}(t)$ is performed using the adaptive Gaussian quadrature routines in \texttt{QuadGK.jl} \cite{quadgk}, while the 3D quadrature required to evaluate the LGF for Mehrstellen stencils is performed with an algorithm from Genz and Malik \cite{genz1980remarks} implemented in \texttt{HCubature.jl} \cite{hcubature}. For domains with one unbounded direction, the expressions given in the appendix are implemented in a C++ header that is available alongside the Julia package. 

\subsection{Residuals in an unbounded domain}\label{subsec:residuals_3unb}
One of the advantages of solving the Poisson equation via convolution with an LGF is that the solution exactly satisfies the discretized system $\mathcal{L}u = \mathcal{R}f$, which may have its own conservation or positivity properties that mimic those of the continuous PDE. To test the degree of exactness in the approximations $\tilde{G}(\vb{n})$ calculated with the algorithms developed in section \ref{sec:poisson}, \ReviewerTwo{we evaluate the residual function 
$R(\vb{n}) \equiv [\mathcal{L} \tilde{G}](\vb{n}) - [\mathcal{R} \delta](\vb{n})$
on a box $[0, N]^3$ and record the maximum absolute residual $\Rmax \equiv \max_{\bm{n} \in [0, N]^3} \abs{R(\bm{n})}$ as well as its location $\vbnres \equiv \argmax_{\bm{n} \in [0, N]^3} \abs{R(\bm{n})}$.}
For dimension split stencils, the region $\abs{\vb{n}} < \nmax$ is precomputed using the algorithm developed in section \ref{sec:unbounded} with an absolute error tolerance of $10^{-15}$, \ReviewerTwo{nearly double precision}. For Mehrstellen stencils, this region is precomputed with 3D adaptive quadrature using a \ReviewerTwo{larger} absolute error tolerance of $10^{-12}$, \ReviewerTwo{which is achievable in a few hours on a laptop computer}. The threshold $\nmax$ is chosen so that a far-field expansion accurate to  $\mathcal{O}({\abs{\vb{n}}^{-19}})$ will meet the same \ReviewerTwo{absolute} error tolerance for $\abs{\vb{n}} > \nmax$, and the remaining $\vb{n}$ that fall in the box $\vb{n}\in [0, 63]^3$ are precomputed using this expansion. All other evaluations of $G(\vb{n})$ are performed with a far-field expansion that is accurate to $\mathcal{O}(\abs{\vb{n}}^{-11})$, which is sufficient for full double precision accuracy. 

Table \ref{tab:3unb_residuals} lists the maximum absolute residual \ReviewerTwo{$\Rmax$} observed within a computational domain of size $N = 128$ \ReviewerTwo{for each stencil considered in this work}, along with the time required to precompute the region $\abs{\vb{n}} \le \nmax$ on a single \ReviewerTwoStrike{core from an Apple M1 Pro processor} \ReviewerTwo{CPU}. \ReviewerTwo{In all cases the location of the maximum residual satisfies $\norm{\vbnres}_2 \approx \nmax$, indicating that the largest residuals occur on the boundary between near-field and far-field evaluation strategies. As expected, the maximum residual for each case is of the same order of magnitude as the absolute error tolerance.} The precomputation times demonstrate a speedup of two orders of magnitude for dimension-split stencils compared to the Mehrstellen stencils, even while the former are computed with significantly tighter error tolerances. \ReviewerOne{While runtimes are of secondary importance for data that can be precomputed and stored, we emphasize that the cost of 3D adaptive quadrature is $\order{N_{\mathrm{sing}}^3}$ where $N_{\mathrm{sing}}$ is the cost of one dimensional quadrature with a singular integrand and specified error tolerance. The algorithm developed in section \ref{sec:three_unbounded_split} reduces this cost to $\order{N^2}$, where $N$ is the cost of quadrature with a smooth integrand and the same error tolerance. This reduction in complexity is essential for computations with tight error tolerances and correspondingly high $N_{\mathrm{sing}}$, and allows for precomputations with near machine-precision accuracy on a laptop computer.}

\begin{table}[htb!]
    \centering
    {\renewcommand{\arraystretch}{1.1}%
    \begin{tabular}{ c  c  c  c  c  c }
        \toprule
        Stencil & $\nmax$ & $N_{\mathrm{evals}}$ & precomputation time & $\Rmax$ & \ReviewerTwo{$\norm{\vbnres}_2$}  \\
        \midrule
        LGF2 & 19 & 779 & $11.5$ s & $2.26\ex{-15}$ &  \ReviewerTwo{19.2} \\ %
        LGF4 & 18 & 672 & $14.1$ s & $2.59\ex{-15}$ &  \ReviewerTwo{18.2} \\ %
        LGF6 & 18 & 672 & $14.5$ s & $2.70\ex{-15}$ &  \ReviewerTwo{18.2} \\ %
        LGF8 & 18 & 672 & $14.9$ s & $2.42\ex{-15}$ &  \ReviewerTwo{18.2}\\ %
        MEH4 & 10 & 141 & $9.43\ex{3}$ s & $1.97\ex{-12}$ & \ReviewerTwo{10.4} \\ %
        MEH6 &  9 & 106 & $6.28\ex{3}$ s & $8.09\ex{-13}$ & \ReviewerTwo{9.0} \\ %
        \bottomrule
    \end{tabular}}
    \caption{Computational cost and maximum residual obtained when evaluating the LGF on a 3D unbounded domain of size $N = 128$. $\nmax$ indicates the radius beyond which the LGF is computed via far-field expansion, and $N_{\mathrm{evals}}$ indicates the number of evaluations necessary to fill the region $\abs{\vb{n}} < \nmax$ after accounting for symmetry. The \ReviewerTwo{final} column indicates the location of the maximum residual within the domain. }
    \label{tab:3unb_residuals}
\end{table}

\subsection{Residuals with one unbounded dimension}\label{subsec:residuals_1unb}
To measure the residuals for LGFs on domains with one unbounded direction, we fix a domain size of $N^3$ points and let $\mathcal{K} = \{2\pi n / N,\, 0 \le n < N\}$ be the set of wavenumbers appearing in a discrete Fourier transform along each dimension of the domain. The partially transformed LGF $G(n_1, k_2, k_3)$ is evaluated for all $0 \le n < N$ and all $k_2,\,k_3 \in \mathcal{K}$ using the analytical expressions provided in Appendix \ref{sec:1unb_expressions}, and the corresponding real-space LGF is defined via the inverse Fourier transform
\begin{equation}\label{eq:transform_1unb}
G(\vb{n}) = \frac{1}{N^2} \sum_{k_2, k_3 \in \mathcal{K}} G(n_1, k_2, k_3) e^{i(n_2k_2 + n_3 k_3)}.
\end{equation}
As discussed in section \ref{sec:one_unbounded}, there is no single analytical expression for $G(n, k_2, k_3)$ that holds uniformly across the domain $k_2, k_3 \in [0, 2\pi]^3$. Particular values of $(k_2, k_3)$ can lead to a variety of singular or pathological behavior, introducing numerical instabilities that must be explicitly addressed. As a result the choice of domain size $N$, and consequently the set of frequencies $\mathcal{K}$ in the Fourier transform, can have a significant effect on the residual. Table \ref{tab:1unb_residuals} lists the maximum absolute residuals encountered on partially unbounded domains of size $N^3$ for $N \in \{30, 56, 176, 416, 768, 1024\}$ and for each stencil listed in section \ref{sec:stencils}. The domains are chosen so that no domain size divides another, leading to corresponding sets $\mathcal{K}$ that contain frequencies unique to that domain. The maximum residual is at most $1.09\ex{-15}$ for dimension-split stencils \ReviewerTwo{across all domain sizes, which is nearly machine precision. For Mehrstellen stencils the maximum residual decreases consistently with increasing domain size, and takes a maximum value of $3.41\ex{-14}$ for MEH4 on a domain of size $30^3$. We attribute the decrease in error to the factor of $1/N^2$ in (\ref{eq:transform_1unb}), which diminishes the influence that any single wavenumber $(k_2, k_3)$ can have on the full solution.}

\begin{table}[htb!]
    \centering
    \begin{tabular}{ c  c  c  c  c  c  c }
        \toprule
        N & LGF2 & LGF4 & LGF6 & LGF8 & MEH4 & MEH6 \\
        \midrule 
        30   & ${3.31\ex{-16}}$ & $8.28\ex{-16}$ & $4.44\ex{-16}$ & $1.09\ex{-15}$ & $3.41\ex{-14}$ & $3.83\ex{-15}$ \\
        56   & $1.38\ex{-16}$ & $4.12\ex{-16}$ & $2.76\ex{-16}$ & $5.04\ex{-16}$ & $1.07\ex{-14}$ & $3.11\ex{-15}$ \\
        176  & $2.22\ex{-16}$ & $1.68\ex{-16}$ & $2.78\ex{-16}$ & $3.31\ex{-16}$ & $3.86\ex{-15}$ & $9.99\ex{-16}$ \\
        416  & $2.22\ex{-16}$ & $2.22\ex{-16}$ & $2.22\ex{-16}$ & $4.44\ex{-16}$ & $1.84\ex{-15}$ & $7.09\ex{-16}$ \\
        768  & $6.25\ex{-17}$ & $2.22\ex{-16}$ & $4.44\ex{-16}$ & $8.88\ex{-16}$ & $1.51\ex{-15}$ & $3.89\ex{-16}$ \\
        1024 & $1.17\ex{-16}$ & $1.46\ex{-16}$ & $4.44\ex{-16}$ & $4.44\ex{-16}$ & $9.65\ex{-16}$ & $4.44\ex{-16}$ \\
        \bottomrule
    \end{tabular}
    \caption{Maximum absolute residuals for domains with one periodic dimension. }
    \label{tab:1unb_residuals}
\end{table}

\subsection{Convergence of Poisson solutions}\label{subsec:convergence}
To demonstrate the utility of LGFs in solving a discretized Poisson equation, we consider the continuous Poisson problem $-\nabla^2 u = f$ with manufactured solution $u = u_1(x_1) u_2(x_2) u_3(x_3)$ and the corresponding right hand side $f = -u''_1(x_1) u_2(x_2) u_3 (x_3) - u_1(x_1) u''_2(x_2) u_3(x_3) - u_1(x_1) u_2(x_2) u_3''(x_3)$. The form of $u_i(x_i)$ varies with the boundary condition along dimension $i$: for periodic or unbounded directions, the corresponding expressions are \ReviewerTwoStrike{the sinusoid and compact Gaussian}
\begin{equation}\label{eq:results_functional_forms}
    u_{per}(x) = \exp(\sin(\frac{8\pi x}{L})) - 1, \quad u_{unb}(x) = \exp(10\qty(1 - \frac{1}{1 - (\frac{2x}{L} - 1)^2})).
\end{equation}
\ReviewerTwo{Similar manufactured solutions are used for validation in \cite{spietz2018regularization, caprace2021flups, balty2023flups}, though we have replaced the sinusoidal $u_{per}(x)$ used by these authors with an exponentiated sinusoid to generate a richer frequency content along each periodic dimension.}
Discretely, we consider a computational domain of size $[0, L]^3$ with grid spacing $h = L / N$. The resulting $N^3$ grid points are placed at coordinates $x_i = (i + \frac{1}{2}) h$ for $0 \le i \le N - 1$ along each axis, and the right hand side $f$ is evaluated at each grid point based on the expressions given in (\ref{eq:results_functional_forms}). The approximate solution is obtained by convolving the right hand side with the LGF for a given finite difference scheme and set of boundary conditions, using an FFT-based convolution provided by the software package FLUPS \cite{caprace2021flups}. The solve is repeated with a range of difference schemes and resolutions, and for each solution we record the maximum pointwise error between the numerical solution and the manufactured solution. 

\begin{figure}
\centering
\begin{subfigure}{0.48\textwidth}
    \centering
    \resizebox{\textwidth}{!}{\input{tikzfigures/three_unbounded.tex}}
    \caption{Fully unbounded domain.}
    \label{fig:convergence_3u}
\end{subfigure}
\begin{subfigure}{0.48\textwidth}
    \centering
    \resizebox{\textwidth}{!}{\input{tikzfigures/one_unbounded}}
    \caption{One dimension unbounded.}
    \label{fig:convergence_1u}
\end{subfigure}
\caption{Convergence data for the 3D Poisson problem with stencils LGF2 (\includegraphics{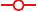}), LGF4 (\includegraphics{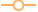}), MEH4 (\includegraphics{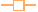}), LGF6(\includegraphics{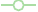}), MEH6 (\includegraphics{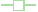}), LGF8 (\includegraphics{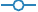}). Dashed lines indicate the expected order of convergence for each difference scheme.}
\label{fig:poisson_results}
\end{figure}

Figure~\ref{fig:poisson_results} provides convergence data for each of the dimension-split and Mehrstellen stencils considered in this work for domains with either one or three unbounded directions and resolutions from $N = 32$ to $N = 1024$. \ReviewerTwo{For the case of three unbounded directions, the LGF is evaluated using the same methodology and absolute error tolerances as in section \ref{subsec:residuals_3unb}}. Each stencil achieves the expected order of convergence on both domain types\ReviewerTwoStrike{, with the Mehrstellen stencils slightly outperforming the dimension-split stencil of equal order at any given resolution}, \ReviewerTwo{with only LGF8 exhibiting a plateau in convergence at higher resolutions. For the fully unbounded domain the plateau occurs at an error of $\epsilon_{\infty} =3.4\ex{-15}$, which is of the same order of magnitude as the absolute error tolerance used in the evaluation of the LGF. For the case of one unbounded direction the plateau begins at $\epsilon_{\infty} = 1.8\ex{-12}$, which is notably larger than the residuals listed for LGF8 in Table \ref{tab:1unb_residuals}. We attribute this discrepancy to an evaluation strategy for $G(n, k_1, k_2)$ that is consistent across all $n$, leading to an LGF evaluation error that varies smoothly in space. This allows a cancellation of error to occurs when the finite difference stencil is applied to the smooth error field to calculate the residual. The error itself we attribute to a loss of precision due to cancellation in the analytical expression given in Appendix \ref{app:lgf8}, which contains polynomials composed of terms with large coefficients and alternating signs. Applications which require accuracy below $10^{-12}$ may require further manipulation of this expression or the use of higher precision arithmetic.}

\section{Conclusions}\label{sec:conclusion}
The algorithms provided in this work extend existing techniques to allow for high-accuracy computations of the LGFs associated with high-order finite difference schemes, and to allow for the calculation of LGFs on domains with one unbounded direction. For dimension-split stencils on unbounded domains, the series expansions and evaluation strategy introduced in section \ref{sec:unbounded} allow for LGF evaluations that are more accurate than 3D quadrature and several hundred times faster. For general stencils on \ReviewerOne{domains with one unbounded direction, we provide numerically stable closed-form expressions} that allow for \ReviewerOne{LGF} evaluations with \ReviewerOne{near} machine-precision accuracy. These algorithms are implemented in an efficient open-source code, and their effectiveness has been demonstrated through numerical experiments.

While we have focused on six stencils in particular, four dimension-split and two Mehrstellen, the provided algorithms are applicable to the \ReviewerTwoStrike{vast }majority of finite difference Poisson discretizations in common use. For unbounded domains, the extension requires minimal effort; the provided open source code automatically generates the necessary asymptotic expansions, compiles them to efficient executable code, and switches between evaluation strategies to \ReviewerOneStrike{minimize the computational cost of precomputing data }\ReviewerOne{achieve a specified error tolerance}. Far field expansions are also automatically generated to any desired precision, and converted to C++ source files. For domains with one unbounded direction the extension is more hands-on, and requires identifying and mitigating numerical instabilities in the analytical expressions for the LGF. However, we believe the techniques introduced to handle the six example stencils are likely to cover the \ReviewerTwoStrike{vast }majority of pathological behavior encountered in practice.

 There are several future directions left open by this work. We have not attempted to \ReviewerOneStrike{accelerate }\ReviewerOne{improve} the numerical quadrature necessary for Mehrstellen stencils on unbounded domains, and there is likely room either for further analytical work or a specialized quadrature algorithm to deliver machine-precision results at a greatly reduced cost. We also have not investigated Mehrstellen stencils beyond sixth order, which have larger stencil widths and could provide characteristic polynomials with pathological behaviors not addressed in section \ref{sec:one_unbounded}. Finally, we have not addressed the computation of Lattice Green's Functions in domains with two unbounded dimensions and one periodic dimension, which appear in the study of vortex tubes and the wakes of cylinders in 3D domains. While the algorithm derived in section \ref{sec:unbounded} can be applied effectively to precompute near-field values, \ReviewerTwo{the derivation of a far-field approximation of the 2D Fourier integral defining the LGF is an open question. Unlike the fully unbounded LGF, the LGF with two unbounded directions is defined by a 2D Fourier integral with a smooth integrand, providing neither a singularity nor a small parameter that could form the basis of such an expansion.  The work of \cite{spietz2018regularization} provides an approximate FFT-based convolution technique for lattices with two unbounded directions that could circumvent this issue; however, as our interest is in exact solutions on unbounded lattices, we have not pursued this approximation further.}

\section*{Acknowledgements}
The authors gratefully acknowledge funding from the Department of Energy Advanced Scientific Computing Research (ASCR) Program, Program Manager Dr. Steven Lee, award number DE-SC0020998.

\appendix
\section{Asymptotic integral expansions for dimension split stencils}\label{section:expand_I_sigma_n_t}
To recover an asymptotic expansion of (\ref{eq:def_I_sigma_n_t}) we rely on the fact that $\sigma(k) - k^2 = \mathcal{O}(k^4)$ as $k \rightarrow 0$. Separating the symbol into its leading order quadratic behavior and higher-order corrections gives
\begin{align}
    I_{\sigma, n}(t) &= \frac{1}{2\pi} \int_{-\pi}^\pi e^{-tk^2} \qty[e^{-t\qty(\sigma(k) - k^2)} \cos(n k)] \dd{k} \nonumber \\
                   &= \frac{1}{2\pi} \int_{-\pi}^\pi e^{-tk^2} \qty[ \sum_{j=0}^\infty \frac{k^{2j}}{(2j)!} a_{2j}(n, t)]\dd{k}
                   \label{eq:no_periodic_gaussian_form}
\end{align}
Here the $a_{2j}(n,t)$ are polynomials in $n$ and $t$, derived by expanding the bracketed quantity as a power series in $k$. Explicitly,
\begin{equation}
    a_{2j}(n,t) = \dv[(2j)]{}{k}\eval_{k=0} \qty[  \exp(-t\qty(\sigma(k) - k^2))\cos(nk) ].
\end{equation}
For large $t$, the Gaussian $e^{-tk^2}$ has support that is concentrated on the origin, so that it is permissible to extend the limits of integration to $[-\infty, \infty]$. With this extension the integral represents a sum over moments of a Gaussian distribution, which can be evaluated in closed form to give
\begin{equation}\label{eq:3unb_series_I_a2j}
    I_{\sigma, n}(t) = \frac{1}{\sqrt{4\pi t}} \sum_{j=0}^\infty \frac{a_{2j}(n, t)}{(2j)!! (2t)^j}.
\end{equation}
The summation above is not a proper power series in $t^{-1}$, since the polynomials $a_{2j}(n, t)$ depend on $t$. After combining like powers of $t$, we arrive at the expansion given in (\ref{eq:3unb_series_I_b2j}) with coefficients $b_{j}(n)$.

In practice, the $b_{j}(n)$ for $j < \jmax$ are computed symbolically with a computer algebra system. The $a_{2j}(n, t)$ for $j < 2\jmax$ are computed by constructing a power series for $-t\qty(\sigma(k) - k^2))$ truncated at $\mathcal{O}(k^{4\jmax})$, exponentiating that series, and multiplying by a power series expansion of $\cos(nk)$ with the same precision. Assembling (\ref{eq:3unb_series_I_a2j}) and consolidating powers of $t$ results in $\jmax$ complete terms of the series in (\ref{eq:3unb_series_I_b2j}), and an additional $\jmax$ partially-formed terms that are discarded.

\section{Maintaining precision at small, nonzero wavenumbers}\label{sec:1unb_small_c}
The expressions for the LGF $G(n,c)$ given in section \ref{sec:one_unbounded_split} become numerically unstable for $c \sim 0$. In this case $p(z;c)$ has a real root $r_1$ that is very near $z = 1$, and the derivative $p'(r_1; c)$ is small. Consequently $p'(r_1; c)$ suffers from catastrophic cancellation when computed in a naive manner. This can be circumvented by writing
\begin{equation}
    p'(r_1; c) = a_w \qty(r_1 - r_1^{-1}) \sum_{i = 2}^{w}(r_1 - r_i)(r_1 - r_i^{-1}),
\end{equation}
which isolates the cancellation in the term $(r_1 - r_1^{-1})$. To calculate this precisely, we note that for small positive values of $c$ there is a root of $q(\lambda) + c = 0$ near $\lambda = 1$, and by repeatedly differentiating the equation $q(\lambda) + c = 0$ the root can be expressed as a power series in $c$. Letting $q_j = q^{(j)}(1)$ for notational simplicity,
\begin{multline}
\lambda - 1 = \frac{c}{2} - \frac{c^2}{16} q_2 + \frac{c^3}{192} \qty(3 q_2^2-2 q_3)+\frac{c^4}{3072}\qty(-4 q_4-15 q_2^3+20 q_3 q_2) \\
+ \frac{c^5}{61440}\qty(-8q_5+40 q_3^2+105 q_2^4+60 q_4 q_2-210 q_3 q_2^2) + \order{c^6}.
\end{multline}
The corresponding root $r$ can then be expressed as $r - 1 = (\lambda - 1) - \sqrt{(\lambda - 1)^2 + 2(\lambda - 1)}$, which does not suffer from catastrophic cancellation. The quantity $(r - r^{-1})$ is then computed precisely via a Taylor expansion about $r = 1$. Similarly, the quantity $r^n$ can also suffer from a loss of precision when $n$ is large and $r \approx 1$. To avoid this we instead compute $\exp(n\log(r))$, obtaining the logarithm with a Taylor expansion about $r = 1$.

\section{Explicit expressions for the LGF with one unbounded dimension}\label{sec:1unb_expressions}

\subsection*{LGF4}
For this kernel, the polynomial $p(z; c)$ has two roots inside the unit circle. The polynomial $q(\lambda) + c = 0$ has two roots given by 
\begin{equation}
    \lambda = 4 \pm \sqrt{9 - 3c}. 
\end{equation}
For $c = 0$, the roots in or on the unit circle are $r = 1$ and $r = 7 - 4\sqrt{3}$, so that the full LGF is
\begin{equation}
    G_4(n) = -\frac{n}{2} + K(r^n - 1), \dbstext{with} r  = 7 - 4\sqrt{3} \dbstext{and} K = \frac{12r^2}{(r - 1)^3(r + 1)}.
\end{equation}
For $c < 10^{-3}$, the root of $q(\lambda) + c = 0$ closest to unity is computed with the power series
\begin{equation}
\lambda = 1+\frac{c}{2}+\frac{c^2}{24}+\frac{c^3}{144}+\frac{5 c^4}{3456}+\frac{7 c^5}{20736}+\frac{7 c^6}{82944}+\frac{11 c^7}{497664}+\order{c^8}.
\end{equation}
At $c^* = 3$, there is a single repeated root $\bar{r} = 4 - \sqrt{15}$ inside the unit circle. For $\abs{c - c^*} < 10^{-5}$ the computation relies on an expansion of $G(n,c)$ about $c = c^*$,
\begin{align*}
\begin{split}
    G^{(0)}(n, c^*) &= \bar{r}^n \qty(\frac{4}{5\sqrt{15}} + \frac{n}{5})
\end{split}\\
\begin{split}
    G^{(1)}(n, c^*) &= -\bar{r}^n \qty(\frac{14}{75\sqrt{15}} + \frac{2}{25} n + \frac{4}{25\sqrt{15}} n^2 + \frac{n^3}{150}) 
\end{split}\\
\begin{split}
    G^{(2)}(n, c^*) &= 2 \bar{r}^n \qty(\frac{901}{18750\sqrt{15}} + \frac{74}{3125} n + \frac{253}{3750\sqrt{15}} n^2 + \frac{23}{3750} n^3 + \frac{1}{250\sqrt{15}} n^4 + \frac{n^5}{15000})
\end{split}\\
\begin{split}
    G^{(3)}(n, c^*) &= -6 \bar{r}^n \left(\frac{37313}{2812500\sqrt{15}} + \frac{11267}{1640625} n + \frac{3167}{140625\sqrt{15}} n^2 + \frac{1}{375} n^3 + \frac{31}{11250\sqrt{15}} n^4 \right. \\
    &\qquad\quad + \left. \frac{31}{281250} n^5 + \frac{1}{28125\sqrt{15}} n^6 + \frac{1}{3150000} n^7 \right)
\end{split}
\end{align*}
For all other $c \in [0,\,c^*]$ there are two real roots, and the LGF is computed using the numerically stable form
\begin{equation}
    G_4(n) = \frac{f(r_1) - f(r_2)}{r_1 - r_2} \dbstext{with} f(z) = \frac{12z^{n+1}}{(z - r_1^{-1}) (z - r_2^{-1})}.
\end{equation}
Finally, for all other $c \in [c^*,\,2\sigmax]$ the two roots form conjugate pair, and the LGF is computed using the oscillatory form given in (\ref{eq:1unb_conjugate_roots}).

\subsection*{LGF6}
For this kernel, the polynomial $p(z;c)$ has one real root and one conjugate pair inside the unit circle. The cubic $q(\lambda) + c = 0$ has one real and two complex roots, given in closed form by 
\begin{equation}
    \lambda = \frac{1}{4} \qty(9 + \xi - \frac{95}{\xi})  \dbstext{with} \xi^3 = 360c - 775 + 60 \sqrt{36c^2 - 155c + 405}.
\end{equation}
Here $\xi$ ranges over all three possible complex values, giving one real root and a conjugate pair. When $c = 0$, the expression reduces to $\lambda = 1$ and $\lambda = \frac{1}{8}(23 \pm 9\sqrt{15} i)$. For small $c$ the root closet to unity is calculated with the power series
\begin{equation}
    \lambda = 1+\frac{c}{2}+\frac{c^2}{24}+\frac{c^3}{720}-\frac{c^4}{1152}-\frac{149 c^5}{518400}-\frac{259 c^6}{6220800}+\frac{163 c^7}{62208000}+\order{c^8}.
\end{equation}
For all other values of $c$ in $[0, 2\sigmax]$ the roots of $p(z;c)$ are well separated. The full LGF is calculated using (\ref{eq:1unb_split}), with the contribution from the conjugate pair reduced to real oscillatory form as in (\ref{eq:1unb_conjugate_roots}).

\subsection*{LGF8}\label{app:lgf8}
For the eight order LGF, the polynomial $q(\lambda) + c = 0$ has four roots given explicitly by
\begin{align*}
    \lambda &= \frac{16}{9} - \frac{\eta}{2} \pm \frac{1}{2}\sqrt{-\frac{434}{81} - \frac{81088}{729\eta} - \eta^2} \dbstext{with} \eta = \pm \sqrt{\frac{\xi}{9} + \frac{3780 c - 4655}{9 \xi} - \frac{434}{81}} \dbstext{and} \\ 
    \xi^3 &=  - 273420 c + 1520225 + 35\sqrt{1968941520 - 879221700 c + 223915104 c^2 - 44089920 c^3}.
\end{align*}
For $c < 10^{-3}$, the root closest to unity is calculated with the power series
\begin{equation}
    \lambda = 1+\frac{c}{2}+\frac{c^2}{24}+\frac{c^3}{720}+\frac{c^4}{40320}+\frac{577 c^5}{3628800}+\frac{389 c^6}{6220800}+\frac{34987 c^7}{3048192000}+\order{c^8}
\end{equation}
The polynomial $p(z;c)$ has a repeated root only for the value $c^* = 3.204471924659898$, which is the single real root of the cubic appearing in the above definition of $\xi$. For $c < c^*$ there are two real roots and one conjugate pair lying inside the unit circle, while for $c^* < c < 2\sigmax$ there are two conjugate pairs inside the unit circle. For $c = c^*$ there is repeated root $r_0$ and a conjugate pair $r_1, \bar{r}_1 = \rho e^{\pm i \theta}$ with
\begin{equation}
    r_0 = 0.1401609439298625, \ \ \rho = 0.12718191072648583
 \dbstext{and} \theta = 1.5912866598657263.
\end{equation}
The expansion of the LGF about $c = c^*$ takes the form
\begin{equation}
    G^{(j)}(n, c^*) = r_0^n \sum_{k = 0}^{2j + 1} a_{jk} n^k + \rho^n \cos(n \theta) \sum_{k = 0}^j u_{jk} n^k  + \rho^n \sin(n\theta) \sum_{k = 0}^j v_{jk} n^k,
\end{equation}
where $a_{jk}$, $u_{jk}$, and $v_{jk}$ are constants derived from (\ref{eq:1unb_repeated_residues}) (see the open source implementation for double-precision values).

\subsection*{Mehrstellen 4}
The numerically stable symbols for the fourth order Mehrstellen stencil are
\begin{align}
    \tilde{\sigma}_{\mathcal{L}}(\bm{y}) &= 4(y_1 + y_2 + y_3) - \frac{8}{3} (y_1 y_2 + y_1 y_3 + y_2 y_3),  \\
    \tilde{\sigma}_{\mathcal{R}}(\bm{y}) &= 1 - \frac{1}{3}(y_1 + y_2 + y_3).
\end{align}
For computations with one unbounded dimension, the coefficients of $p_{\mathcal{L}}(z; \bm{y})$ are
\begin{equation}
    a_0(\bm{y}) = 2 - \frac{8}{3}\qty(y_2 y_3 - y_2 - y_3), \quad
    a_1(\bm{y}) = -1 + \frac{2}{3} \qty(y_2 + y_3),
\end{equation}
while the coefficients of $p_{\mathcal{R}}(z; \bm{y})$ are
\begin{equation}
    b_0(\bm{y}) = \frac{5}{6} - \frac{1}{3}\qty(y_2 + y_3), \quad
    b_1(\bm{y}) = \frac{1}{12}.
\end{equation}
The LGF is given in closed form by (\ref{eq:1unb_general}), and after rearranging for numerical stability
\begin{equation}
    G(n, \bm{0}) = -\frac{n}{2} - \frac{1}{12}\delta(n), \qquad 
    G(n, \bm{y}) = \begin{dcases}
         \frac{r^{n} p_{\mathcal{R}}(r;\bm{y})}{a_1(\bm{y}) (r^2 - 1)}, & n \neq 0, \\
        \frac{r p_{\mathcal{R}}'(r;\bm{y})}{a_1(\bm{y}) (r^2 - 1)}, & n = 0.
    \end{dcases}
\end{equation}
Here $r$ is the root of $p_{\mathcal{L}}(z; \bm{y})$ corresponding to $\lambda = -a_0(\bm{y}) / 2a_1(\bm{y})$. When $y_1 + y_2 = \frac{3}{2}$ the polynomial $p_{\mathcal{L}}(z; \bm{y})$ has a root at $z = 0$, and the above should be replaced by
\begin{equation}
    G(n, \vb{y}) = \delta(n) \frac{b_0(\bm{y})}{a_0(\bm{y})} + \delta(n - 1) \frac{b_1(\bm{y})}{a_0(\bm{y})}.
\end{equation}

\subsection*{Mehrstellen 6}
The numerically stable symbols for the sixth order Mehrstellen stencil are
\begin{align}
    \tilde{\sigma}_{\mathcal{L}}(\bm{y}) &=  4(y_1 + y_2 + y_3) -\frac{8}{3} (y_1 y_2 + y_1 y_3 + y_2 y_3) + \frac{32}{15} y_1 y_2 y_3, \\
    \tilde{\sigma}_{\mathcal{R}}(\bm{y}) &= 1 - \frac{1}{3}(y_1 + y_2 + y_3) - \frac{1}{15}(y_1^2 + y_2^2 + y_3^2) + \frac{8}{45}(y_1 y_2 + y_1 y_3 + y_2 y_3).
\end{align}
For computations with one unbounded dimension, the coefficients of $p_{\mathcal{L}}(z; \bm{y})$ are
\begin{equation}
    a_0(\bm{y}) = 2 + \frac{8}{3}\qty(y_2 + y_3) - \frac{8}{5}y_2 y_3 , \quad
    a_1(\bm{y}) = -1 + \frac{2}{3} \qty(y_2 + y_3) - \frac{8}{15} y_2 y_3,
\end{equation}
while the coefficients of $p_{\mathcal{R}}(z; \bm{y})$ are
\begin{equation}
    b_0(\bm{y}) =  -\frac{1}{15}(y_2^2 + y_3^2) + \frac{8}{45} y_2 y_3 - \frac{11}{45} (y_2 + y_3) + \frac{97}{120}, \quad
    b_1(\bm{y}) = -\frac{2}{45} (y_2 + y_3) + \frac{1}{10}, \quad
    b_2(\bm{y}) = -\frac{1}{240}.
\end{equation}
The LGF in closed from is then given by
\begin{equation}
        G(n, \bm{y}) =  \frac{r^{n - 1} p_{\mathcal{R}}(r; \bm{y})}{a_1(\bm{y}) (r^2 - 1)} + \delta(n) \frac{b_1(\bm{y}) a_1(\bm{y}) - b_2(\bm{y}) a_0(\bm{y})}{a_1(\bm{y})^2} + \delta(n - 1)\frac{b_2(\bm{y})}{a_1(\bm{y})},
\end{equation}
where $r$ is the root of $p_{\mathcal{L}}(z; \bm{y})$ corresponding to $\lambda = -a_0(\bm{y}) / 2a_1(\bm{y})$.

\bibliography{references}

\begin{thebibliography}{10}

\bibitem{balty2023flups}
{\sc P.~Balty, P.~Chatelain, and T.~Gillis}, {\em Flups-a flexible and
  performant massively parallel fourier transform library}, IEEE Transactions
  on Parallel and Distributed Systems,  (2023).

\bibitem{beckers2022planar}
{\sc D.~Beckers and J.~D. Eldredge}, {\em Planar potential flow on cartesian
  grids}, Journal of Fluid Mechanics, 941 (2022), p.~A19.

\bibitem{buneman1971analytic}
{\sc O.~Buneman}, {\em Analytic inversion of the five-point poisson operator},
  Journal of Computational Physics, 8 (1971), pp.~500--505.

\bibitem{caprace2021flups}
{\sc D.-G. Caprace, T.~Gillis, and P.~Chatelain}, {\em Flups: A fourier-based
  library of unbounded poisson solvers}, SIAM Journal on Scientific Computing,
  43 (2021), pp.~C31--C60.

\bibitem{chatelain2013}
{\sc P.~Chatelain, S.~Backaert, G.~Winckelmans, and S.~Kern}, {\em Large eddy
  simulation of wind turbine wakes}, Flow, turbulence and combustion, 91
  (2013), pp.~587--605.

\bibitem{chatelain2010fourier}
{\sc P.~Chatelain and P.~Koumoutsakos}, {\em A fourier-based elliptic solver
  for vortical flows with periodic and unbounded directions}, Journal of
  Computational Physics, 229 (2010), pp.~2425--2431.

\bibitem{cserti2000application}
{\sc J.~Cserti}, {\em Application of the lattice green’s function for
  calculating the resistance of an infinite network of resistors}, American
  Journal of Physics, 68 (2000), pp.~896--906.

\bibitem{deriaz2020compact}
{\sc E.~Deriaz}, {\em Compact finite difference schemes of arbitrary order for
  the poisson equation in arbitrary dimensions}, BIT Numerical Mathematics, 60
  (2020), pp.~199--233.

\bibitem{dorschner2020fast}
{\sc B.~Dorschner, K.~Yu, G.~Mengaldo, and T.~Colonius}, {\em A fast
  multi-resolution lattice green's function method for elliptic difference
  equations}, Journal of Computational Physics, 407 (2020), p.~109270.

\bibitem{eldredge2022method}
{\sc J.~D. Eldredge}, {\em A method of immersed layers on cartesian grids, with
  application to incompressible flows}, Journal of Computational Physics, 448
  (2022), p.~110716.

\bibitem{AbstractAlgebra.jl-2017}
{\sc C.~Fieker, W.~Hart, T.~Hofmann, and F.~Johansson}, {\em Nemo/hecke:
  Computer algebra and number theory packages for the julia programming
  language}, in Proceedings of the 2017 ACM on International Symposium on
  Symbolic and Algebraic Computation, ISSAC '17, New York, NY, USA, 2017, ACM,
  pp.~157--164, \url{https://doi.org/10.1145/3087604.3087611},
  \url{https://doi.acm.org/10.1145/3087604.3087611}.

\bibitem{gabbard2022immersed}
{\sc J.~Gabbard, T.~Gillis, P.~Chatelain, and W.~M. van Rees}, {\em An immersed
  interface method for the 2d vorticity-velocity navier-stokes equations with
  multiple bodies}, Journal of Computational Physics,  (2022), p.~111339.

\bibitem{genz1980remarks}
{\sc A.~C. Genz and A.~A. Malik}, {\em Remarks on algorithm 006: An adaptive
  algorithm for numerical integration over an n-dimensional rectangular
  region}, Journal of Computational and Applied mathematics, 6 (1980),
  pp.~295--302.

\bibitem{gholami2016fft}
{\sc A.~Gholami, D.~Malhotra, H.~Sundar, and G.~Biros}, {\em Fft, fmm, or
  multigrid? a comparative study of state-of-the-art poisson solvers for
  uniform and nonuniform grids in the unit cube}, SIAM Journal on Scientific
  Computing, 38 (2016), pp.~C280--C306.

\bibitem{gillis2018fast}
{\sc T.~Gillis, G.~Winckelmans, and P.~Chatelain}, {\em Fast immersed interface
  poisson solver for 3d unbounded problems around arbitrary geometries},
  Journal of Computational Physics, 354 (2018), pp.~403--416.

\bibitem{gillman2010fast}
{\sc A.~Gillman and P.-G. Martinsson}, {\em Fast and accurate numerical methods
  for solving elliptic difference equations defined on lattices}, Journal of
  Computational Physics, 229 (2010), pp.~9026--9041.

\bibitem{gowda2022}
{\sc S.~Gowda, Y.~Ma, A.~Cheli, M.~Gw\'{o}\'{z}zd\'{z}, V.~B. Shah, A.~Edelman,
  and C.~Rackauckas}, {\em High-performance symbolic-numerics via multiple
  dispatch}, ACM Commun. Comput. Algebra, 55 (2022), p.~92–96,
  \url{https://doi.org/10.1145/3511528.3511535},
  \url{https://doi.org/10.1145/3511528.3511535}.

\bibitem{hejlesen2013high}
{\sc M.~M. Hejlesen, J.~T. Rasmussen, P.~Chatelain, and J.~H. Walther}, {\em A
  high order solver for the unbounded poisson equation}, Journal of
  Computational Physics, 252 (2013), pp.~458--467.

\bibitem{hejlesen2015high}
{\sc M.~M. Hejlesen, J.~T. Rasmussen, P.~Chatelain, and J.~H. Walther}, {\em
  High order poisson solver for unbounded flows}, Procedia IUTAM, 18 (2015),
  pp.~56--65.

\bibitem{hockney1988computer}
{\sc R.~Hockney and J.~Eastwood}, {\em Computer simulation using particles
  taylor \& francis}, Inc., USA,  (1988).

\bibitem{ji2022bursting}
{\sc L.~Ji and W.~M. Van~Rees}, {\em Bursting on a vortex tube with initial
  axial core-size perturbations}, Physical Review Fluids, 7 (2022), p.~044704.

\bibitem{hcubature}
{\sc S.~G. Johnson}, {\em {HCubature.jl}}.
\newblock \url{https://github.com/JuliaMath/HCubature.jl}.

\bibitem{quadgk}
{\sc S.~G. Johnson}, {\em {QuadGK.jl}: {G}auss--{K}ronrod integration in
  {J}ulia}.
\newblock \url{https://github.com/JuliaMath/QuadGK.jl}, 2013.

\bibitem{kavouklis2018computation}
{\sc C.~Kavouklis and P.~Colella}, {\em Computation of volume potentials on
  structured grids with the method of local corrections}, Communications in
  Applied Mathematics and Computational Science, 14 (2018), pp.~1--32.

\bibitem{koster1954simplified}
{\sc G.~Koster and J.~Slater}, {\em Simplified impurity calculation}, Physical
  Review, 96 (1954), p.~1208.

\bibitem{liska2014parallel}
{\sc S.~Liska and T.~Colonius}, {\em A parallel fast multipole method for
  elliptic difference equations}, Journal of Computational Physics, 278 (2014),
  pp.~76--91.

\bibitem{liska2016fast}
{\sc S.~Liska and T.~Colonius}, {\em A fast lattice green's function method for
  solving viscous incompressible flows on unbounded domains}, Journal of
  Computational Physics, 316 (2016), pp.~360--384.

\bibitem{martinsson2002asymptotic}
{\sc P.-G. Martinsson and G.~J. Rodin}, {\em Asymptotic expansions of lattice
  green's functions}, Proceedings of the Royal Society of London. Series A:
  Mathematical, Physical and Engineering Sciences, 458 (2002), pp.~2609--2622.

\bibitem{spietz2018regularization}
{\sc H.~J. Spietz, M.~M. Hejlesen, and J.~H. Walther}, {\em A regularization
  method for solving the poisson equation for mixed unbounded-periodic
  domains}, Journal of Computational Physics, 356 (2018), pp.~439--447.

\bibitem{zucker201170+}
{\sc I.~Zucker}, {\em 70+ years of the {Watson} integrals}, Journal of
  Statistical Physics, 145 (2011), pp.~591--612.

\end{thebibliography}

\end{document}